\documentclass[12pt,psamsfonts]{amsart} \usepackage{amssymb}
\usepackage{amscd} \usepackage{longtable}
\calclayout

\numberwithin{equation}{subsection}

\newtheorem{theorem}[subsection]{Theorem}
\newtheorem{lemma}[subsection]{Lemma}
\newtheorem{proposition}[subsection]{Proposition}
\newtheorem{corollary}[subsection]{Corollary}

\theoremstyle{definition}
\newtheorem{definition}[subsection]{Definition}

\theoremstyle{remark} \newtheorem{remark}[subsection]{Remark}

\newtheorem{example}[subsection]{Example}

\def\inv{^{-1}} \def\Ker{\operatorname{Ker}}
 \def\var{\varphi}
 
 \def\today{\ifcase\month\or January\or
  February\or March\or April\or May\or June\or July\or August\or
  September\or October\or November\or December\fi \space\number\day,
  \number\year} \def\C{\mathbb C} 
\def\quot#1#2{{#1/\!\!/#2}} \def\twist#1#2#3{#1\times^{#2}#3}
\def\c{^\C} \def\xiz{{\xi_Z}} \def\lie#1{\mathfrak{ #1}}
 \def\liek{\lie k} \def\liep{\lie p} \def\lieq{\lie
  q} \def\liea{\lie a} \def\lieg{\lie g} \def\liem{\lie m}
\def\lieu{\lie u}  
\def\R{{\mathbb R}} \def\Z{{\mathbb Z}}
\newcommand{\Ad}{\operatorname{Ad}} \def\muip{\mu_{i\liep}}
\def\muk{\mu_{\liek}} \def\mip{\mathcal M_{i\liep}} \def\m{\mathcal M}
 \def\tosim{\xrightarrow{\sim}} 
\def\pt{\partial} \def\Comp{\operatorname{Comp}}
\def\SU{\operatorname{SU}} \def\SL{\operatorname{SL}}
\def\tr{\operatorname{tr}} 
\def\Ann{\operatorname{Ann}} \def\GL{\operatorname{GL}}
\def\Hom{\operatorname{Hom}} \def\Sing{\operatorname{Sing}}
\title{Cartan decomposition of the moment map} \author{Peter Heinzner}
\thanks{First author is partially supported by the
  Sonderforschungsbereich SFB/TR12 of the Deutsche
  Forschungsgemeinschaft}
\address{Fakult\"at f\"ur Mathematik\\
  Ruhr Universit\"at Bochum\\
  Universit\"atsstrasse 150\\
  D - 44780 Bochum} \email{heinzner@cplx.rub.de} \author{Gerald W.
  Schwarz} \thanks{Second author is partially supported by NSA grant
  H98230--04--01--0070}
\address{Department of Mathematics\\
  Brandeis University\\
  PO Box 549110\\
  Waltham, MA 02454-9110} \email{schwarz@brandeis.edu}

\begin{document}
\begin{abstract}  
  We investigate a class of actions of real Lie groups on complex
  spaces.  Using moment map techniques we establish the existence of a
  quotient and a version of Luna's slice theorem as well as a version
  of the Hilbert-Mumford criterion.  A global slice theorem is proved
  for proper actions.  We give new proofs of results of Mostow on
  decompositions of groups and homogeneous spaces.
\end{abstract}
\maketitle

\date{today}
\section{Introduction}

Let $Z$ be a complex space with a holomorphic action of the complex
reductive group $U\c$, where $U\c$ is the complexification of the
compact Lie group $U$.  We assume that $Z$ admits a smooth
$U$-invariant K\"ahler structure and a $U$-equivariant moment mapping
$\mu\colon Z\to\lieu^*$, where $\lieu$ is the Lie algebra of $U$ and
$\lieu^*$ its dual.  We assume that $G\subset U\c$ is a closed
subgroup such that the Cartan decomposition $U\c=U \exp(i\lieu)\simeq
U\times i\lieu$ induces a Cartan decomposition $G=K\exp\liep\simeq
K\times \liep$ where $K=U\cap G$ and $\liep\subset i\lieu$ is an $(\Ad
K)$-stable linear subspace.  We have the subspace $i\liep\subset\lieu$
and by restriction an induced ``moment'' mapping $\muip\colon
Z\to(i\liep)^*$.  We define $\mip$ to be the set of zeroes of $\muip$,
and we define $\m$ to be the set of zeroes of $\mu$.  For a given
$\mu$ we have the set $\mathcal{S}_{U\c}(\m):=\{z\in Z;\ 
\overline{U\c\cdot z}\cap \m\neq\emptyset\}$ of semistable points with
respect to $\mu$ and the $U\c$-action on $Z$. We call
$\mathcal{S}_G(\mip):=\{z\in Z;\ \overline{G\cdot z}\cap
\mip\neq\emptyset\}$ the set of semistable points of $Z$ with respect
to $\muip$ and the $G$-action on $Z$. 

The set $\m$ plays an important role in determining the 
closed $U\c$-orbits in $\mathcal{S}_{U\c}(\m)$ and its quotient under the $U\c$-action 
(see Theorems~\ref{Uc quotients} and~\ref{holomorphic slice theorem}
below). We show that $\mip$ is the right analogue
of $\m$ for the action of $G$, as follows.

\begin{theorem} \label{main introd theorem} 
  Let $Z$, $G$, $\mip$, $\m$ be as above and set $Z':=\mathcal{S}_G(\mip)$. 
\begin{enumerate}
\item \begin{enumerate}
\item (Corollary~\ref{closed orbits of G}) An orbit $G\cdot z$ is
  closed in $Z'$ if and only if $G\cdot z\cap\mip\neq\emptyset$.
\item (Theorem~\ref{first main theorem}) There is a quotient space
  $\quot {Z'}G$ which parametrizes the closed $G$-orbits in $Z'$. The
  inclusion $\mip\to Z'$ induces a homeomorphism $\mip/K\simeq \quot
  {Z'}G$.
\item (Corollary~\ref{minimal dimension orbit}) Let $z\in Z'$. Then 
  $\overline{G\cdot z}\cap Z'$ 
contains a unique orbit of minimal dimension
  and this orbit is closed in $Z'$.
\end{enumerate}

\item (Lemma~\ref{isotropy decomposition G}) Let $z\in\mip$. Then
  $G\cdot z\cap\mip=K\cdot z$ and $G_z=K_z\cdot\exp(\liep_z)$ where
  $\liep_z$ denotes the elements of $\liep$ such that the
  corresponding vector field on $Z$ vanishes at $z$.
\item (Theorem~\ref{slice theorem}, Remark~\ref{saturated slice
    remark}) Let $z\in \mip$. Then there is a locally closed real
  analytic $G_z$-stable subset $S$ of $Z$, $z\in S$, such that the
  natural map $\twist G{G_z}S\to Z$ is a real analytic $G$-isomorphism
  onto the open set $G\cdot S$. Moreover, if $Z'$ 
  is open, then $S$ can be chosen such that $G\cdot
  S$ is saturated with respect to the quotient map $Z'\to \quot {Z'}G$.
\item If $Z=\mathcal{S}_{U\c}(\m)$, then
\begin{enumerate}
\item (Proposition~\ref{G-orbit closure}) $Z=\mathcal{S}_G(\mip)$.
\item (Corollary~\ref{hilbert for analytic quotient}) Let $z\in Z$ and
  suppose that $Y\subset \overline{G\cdot z}$ is closed and
  $G$-stable. Then there is a Lie group homomorphism
  $\lambda\colon\R\to G$ such that $\lim_{t\to -\infty}
  \lambda(t)\cdot z$ exists and is a point in $Y$.  The image of
  $\lambda$ consists of semisimple elements of $G$.
\end{enumerate}

\end{enumerate}
\end{theorem}

Note that (1b) implies that   $\quot {Z'}G$ is  
metrizable and locally compact.
If $Z'$ is open, then from (3) one can deduce that  $\quot {Z'}G$ is 
locally homeomorphic to
semianalytic sets (Corollary~\ref{semianalytic corollary}).

If $Z$ is
a Stein space then it admits a smooth strictly plurisubharmonic
$U$-invariant exhaustion function $\rho$. Associated with $\rho$ is a
$U$-invariant K\"ahler structure and a moment mapping $\mu$. 
Then $\m$ is the Kempf-Ness set (see
\cite{KempfNess, Schwarz}). 
Moreover, for any such $\mu$, the equality $Z=\mathcal{S}_{U\c}(\m)$
holds automatically. Another interesting example of equality is the
case where $Z$ is the set of semistable points (in the sense of
geometric invariant theory) relative to a $U\c$-linearized ample line
bundle of a projective variety $Z_0$.  In this case there also exists a
$U$-invariant K\"ahler structure on $Z$ and a $\mu$ such that
$Z=\mathcal{S}_{U\c}(\m)$. Moreover, $Z_0':=\mathcal{S}_G(\mip)$ is then a 
$G$-stable open subset of $Z_0$
which is not usually $U\c$-stable (Remark \ref{sgmip open remark}).

Of course, there is much earlier work on quotients and slice theorems
for actions of complex reductive groups, and there is also earlier
work for actions of real groups. In particular, in the latter case,
there are the papers of Richardson-Slodowy~\cite{RichardsonSlodowy}
and Luna~\cite{LunaReal}. Here one has a complex representation space
$V$ of $U\c$ and real forms $V_\R$ of $V$ and $G$ of $U\c$. One
considers the action of $G$ on $V_\R$. Here $V_\R$ is a Lagrangian
subspace of $V$.  In the case where $G$ is a real form of $U\c$ and
$Z$ is a K\"ahler manifold there are also results about the structure
of the $G$-action on Lagrangian submanifolds $X$ of $Z$ using moment
map techniques (see, e.g., \cite{OsheaSjamaar} and references
therein).  These cases are rather special.  The $\muip$-component of
$\mu$ on $X$ is completely determined by $\mu$.  One establishes
results concerning $\m$ and the $U\c$-action on $Z$ and then restricts
to $X$. This works because $X\cap \m=X\cap \mip$ and because the map
$\muk\colon Z\to \liek^*$ obtained by restricting $\mu$ to $\liek$ is
constant on $X$.

The results presented here are much more general: the group actions 
are not necessarily algebraic, the group $G$ is not necessarily a 
real form of $U\c$ and
we consider the action on $Z$, not just on a real form of $Z$.

Besides the results mentioned above, we also consider several topics
pertaining to proper actions and compact isotropy groups.  In
particular we show the following.

\begin{theorem} \label{introd theorem}
  (Proposition~\ref{real proper}, Remark~\ref{real proper remark})
  Assume that $Z=\mathcal{S}_{U\c}(\m)$. Let $X$ be a $G$-stable
  closed subset of $Z$ such that the $G$-action on $X$ is proper. Then
  the natural map $\twist GK(\mip\cap X)\to X$ is a homeomorphism and
  a real analytic isomorphism if $X$ and $\muip$ are real analytic.
\end{theorem}

We have a similar decomposition for the subset $\Comp_{i\liep}(Z)$ of
points $z\in Z$ such that $G_z$ is compact (Theorem~\ref{compact slice
  theorem}).  The results on proper actions are applied to obtain
decompositions, due to Mostow, for groups and homogeneous spaces. The
application relies on properties of a distinguished strictly
plurisubharmonic exhaustion of $U\c$ related to the Cartan
decomposition. See section~\ref{decomps of homogeneous spaces}.

Several of our results rely upon the notion of $\muip$-adapted sets.
A $\muip$-adapted subset of $Z$ is a $K$-invariant subset $A$ of $Z$
such that for all $z\in Z$ and $\xi\in i\liep$, the curve $(\exp
it\xi)\cdot z$ lies in $A$ for a connected set $J$ of $t\in \R$.
Moreover, we require that if $t_+:=\sup J<\infty$, then
$\muip(\exp(it_+\xi)\cdot z)(i\xi)>0$ and a similar negativity
condition if $t_-:=\inf J>-\infty$. 
We are able to show that every $K$-orbit in $\mip$ has a neighborhood basis of open $\muip$-adapted sets in the case
that $Z=\mathcal{S}_{U\c}(\m)$ (Theorem~\ref{main theorem}). The
$\muip$-adapted sets have very nice properties. For example, if $A_1$
and $A_2$ are $\muip$-adapted, then $G\cdot A_1\cap G\cdot
A_2=G\cdot(A_1\cap A_2)$.

In case $U$ is commutative, 
we can prove Theorem~\ref{main theorem} without the hypothesis $Z=\mathcal{S}_{U\c}(\m)$.
This allows us to establish the
``separation property''~\ref{defn separation property} which
is used in the proof of most of the statements in Theorem~\ref{main introd
  theorem}. In
addition, if $U$ is commutative, then $\mathcal{S}_G(\mip)$ is open 
and $U\c$-stable in $Z$. In general $\mathcal{S}_G(\mip)$ is 
not $U\c$-stable. It would be interesting to know if 
$\mathcal{S}_G(\mip)$ is always open in $Z$.  
Examples indicate that it should be extremely
interesting to clarify the interplay of the various geometric objects
associated with $\mu$, $\muip$ and $\muk$.

The authors thank Henrik St\"otzel for helpful discussions and
remarks.  The second author is thankful for the hospitality of the
Mathematics Department of the Ruhr Universit\"at Bochum where much of
this research took place. 
\section{$G$-fiber bundles and slices}

Let $H$ be a closed subgroup of the Lie group $G$ and $S$ an
$H$-space. The \emph{twisted product} $\twist{G}H{S}$ is the quotient
$(G\times S)/H$ with respect to the $H$-action $(h,(g,x))\mapsto
(gh^{-1}, h\cdot x)$. Since the $G$-action on $G\times S$ given by
multiplication from the left on the first factor commutes with the
$H$-action there is an induced $G$-action on $\twist{G}{H}{S}$. We use
the notation $[g,x]$ for $H\cdot (g,x)\in \twist{G}{H}{S}$.  Note that
$\twist{G}H{S}$ is a $G$-fiber bundle over $G/H$ with fiber $S$
associated to the $H$-principal bundle $G\to G/H$.  Let $X$ be a
$G$-space and $H$ a closed subgroup of $G$. An $H$-stable subspace $S$
of $X$ is said to be a {\it global $H$-slice\/} if the natural map
$\twist{G}{H}{S}\to X$, $[g,z]\mapsto g \cdot z$ is an isomorphism. If
$x\in S$, 
$S$ is $G_x$-stable, $G\cdot S$ is open in $X$  
and $\twist{G}{G_x}{S}\to
G\cdot S$ is an isomorphism, then $S$ is called a \emph{geometric
  slice at $x$}.

\section{K\"ahler structures}

In this paper a complex space $Z$ is always a reduced complex space
with countable topology. If $G$ is a Lie group, then a \emph{complex
  $G$-space $Z$\/} is a complex space with a real analytic action
$G\times Z\to Z$ which for fixed $g\in G$ is holomorphic. For a
complex Lie group $G$ a \emph{holomorphic $G$-space\/} is a complex
$G$-space $Z$ such that the $G$-action $G\times Z\to Z$ is
holomorphic.

A K\"ahler structure $\omega$ on $Z$ is an open covering
$\{U_\alpha\}$ of $Z$ together with smooth strictly plurisubharmonic
functions $\rho_\alpha:U_\alpha\to \R$ such that
$h_{\alpha\beta}=\rho_\alpha-\rho_\beta$ is pluriharmonic on
$U_\alpha\cap U_\beta$. Here strictly plurisubharmonic means strictly
plurisubharmonic with respect to perturbations (see
\cite{HeinznerHuckleberryLoose}) and a pluriharmonic function is, as
usual, a function which is locally the real part of a holomorphic
function. Note that for smooth $Z$ one obtains the usual definition of
a K\"ahler manifold whose K\"ahler form is given locally by
$\omega_\alpha=-dd^c \rho_\alpha=2i\pt\overline{\pt}\rho_\alpha$.
Here $\pt$ and $\overline{\pt}$ are the usual exterior differential
operators and $d^c\rho(v)=d\rho(Jv)$ for every $v\in T_zZ$ and smooth
function $\rho$, where $J$ denotes multiplication by $i=\sqrt{-1}$ on
the tangent space $T_zZ$. For smooth $Z$ we will not distinguish
between $\omega=\{\rho_\alpha\}$ as defined above and the associated
K\"ahler form $\omega$ given by $\omega|{U_\alpha}:=-dd^c\rho_\alpha$.

For a complex $G$-space $Z$ one has a natural notion of a
$G$-invariant K\"ahler structure $\omega$. A moment map on a complex
$G$-space $Z$ with respect to such an invariant K\"ahler structure is
a smooth $G$-equivariant map $\mu \colon Z\to \lieg^*$ such that, for
every $G$-stable complex submanifold $Y$ of $Z$ and $\xi\in \lieg$, we
have

$$d\mu^\xi=\imath_{\xi_Z}\omega_Y.$$

\noindent  Here $\omega_Y$ denotes the K\"ahler form induced on
$Y$ and $\xi_Z$ denotes the vector field on $Z$ induced by $\xi$ and
the $G$-action.  The map $\mu^\xi\colon Z\to\C$ sends $z\in Z$ to
$\mu(z)(\xi)$ and $\imath_{\xi_Z}$ denotes contraction with $\xi_Z$.
Our condition on $d\mu^\xi$ is equivalent to requiring that
$\mathrm{grad} (\mu^\xi)=J\xi_Z$ where the gradient is with respect to
the underlying Riemannian structure on $Y$.

Let $G$ be a complex reductive group, let $Z$ be a holomorphic
$G$-space and let $\mathcal{O}_Z$ denote the structure sheaf of $Z$. A
complex space $Y$ together with a holomorphic map $\pi\colon Z\to Y$
is said to be an \emph{analytic Hilbert quotient of $Z$\/} with
respect to the $G$-action if

\begin{itemize}
\item $\pi$ is a $G$-invariant locally Stein map and
  
\item $\mathcal{O}_Y=\pi_*\mathcal{O}_Z^G$.
\end {itemize}
Here locally Stein means that there is an open covering of $Y$ by open
Stein subspaces $Q_\alpha$ such that $\pi\inv(Q_\alpha)$ is a Stein
subspace of $Z$ for all $\alpha$, and $\pi_*\mathcal{O}^G_Z$ denotes
the sheaf $Q\to \mathcal{O}_Z(\pi\inv(Q))^G$, $Q$ open in $Y$.  If the
analytic Hilbert quotient of $Z$ with respect to $G$ exists it will be
denoted by $\quot{Z}{G}$.  If $H$ is a reductive algebraic subgroup of
$G$ and $\quot ZG$ exists, then $\quot ZH$ exists (see, e.g.,
\cite{HeinznerMiglioriniPolito}).

Assume now that the analytic Hilbert quotient $\quot ZG$ exists. Then
it has the following properties.

\begin{itemize}
\item For every Stein subspace $A$ of $\quot{Z}{G}$ the inverse image
  $\pi\inv(A)$ is a Stein subspace of $Z$.
  
\item For every closed analytic $G$-invariant subspace $X$ of $Z$ the
  image $\pi(X)$ is a closed analytic subspace of $\quot{Z}{G}$ and
  the restriction $\pi|X\colon X\to \pi(X)$ is an analytic Hilbert
  quotient of $X$.
  
\item The quotient map $\pi$ maps disjoint closed $G$-stable subsets
  of $Z$ onto disjoint closed subsets of $\quot{Z}{G}$.

\end{itemize}

If $Z$ is a Stein space, then the analytic Hilbert quotient exists and
has the properties above (see \cite{Heinzner}).  See
\cite{HeinznerMiglioriniPolito} for the general case.

\section{The Cartan decomposition} \label{Cartan decomposition}

Let $U$ be a compact Lie group. Then $U$ has a natural real linear algebraic group structure, 
and we denote by
$U\c$ the corresponding complex linear algebraic group \cite{Chevalley}. The group $U\c$ is
reductive and is the universal complexification of $U$ in the sense of
\cite{Hochschild}.  On the Lie algebra level we have the Cartan
decomposition
$$\lieu\c=\lieu +i\lieu$$
with a corresponding Cartan involution
$\theta:\lieu\c\to\lieu\c$, $\xi+i\eta\mapsto \xi - i\eta$, $\xi$,
$\eta\in\lieu$. We also use $\theta$ to denote the corresponding
anti-holomorphic involution on $U\c$. The real analytic map
$$U\times i\lieu\to U\c, \qquad (u,\xi)\mapsto u\exp \xi$$
is a
diffeomorphism. Since
$g_1(u\exp\xi)g_2\inv=g_1ug_2\inv\exp(\mathrm{Ad}(g_2)\cdot \xi)$ for
$(g_1,g_2)\in U\times U$, $\xi\in\lie u$ and $u\in U$, the isomorphism
$U\times i\lieu\cong U\c$ is a $U\times U$-equivariant diffeomorphism.
We refer to the decomposition $U\c=U\exp(i\lieu)$ as the \emph{Cartan
  decomposition of $U\c$}, and we fix it for the remainder of this
paper.
\\

Let $G$ be a real Lie subgroup of $U\c$. We say that $G$ is
\emph{compatible with the Cartan decomposition of} $U\c$ if $G=K\exp
\liep$ where $K$ is a Lie subgroup of $U$ and $\liep$ is a $K$-stable
linear subspace of $i\lieu$. Note that
$$K\times \liep\to G, \qquad (k, \xi)\to k\exp\xi$$
is a $(K \times
K)$-equivariant diffeomorphism. The closure $\bar G$ of $G$ in $U\c$
is given by $\bar K\exp\liep$. In particular, $G$ is a closed subgroup
of $U\c$ if and only if $K$ is compact. For a closed subgroup $G$ the
complexification $K\c$ of $K$ is a closed complex subgroup of $U\c$
and
$$K\times i\liek\to K\c, \qquad (k, \xi)\to k\exp\xi$$
is a $(K \times
K)$-equivariant diffeomorphism.

\begin{example}
\begin{itemize}
\item[a)] For any compact subgroup $K$ of $U$, both $K$ and its
  complexification $G=K\c$ are compatible with the Cartan
  decomposition of $U\c$. In particular, $G=U\c$ is an example of a
  compatible subgroup.
\item[b)] For any $\xi\in i\lieu$, the group $G=\exp(\R\xi)$ is
  compatible. More generally, if $\liea\subset i\lieu$ is a Lie
  subalgebra of $\lieu\c$, then it is commutative and
  $G=\exp(i\liea+\liea)$ is a compatible subgroup of $U\c$. Note that
  $K=\exp(i\liea)$ need not be compact.
\item[c)] 
Let $\sigma$ be an antiholomorphic
  involution of $U\c$ which commutes with $\theta$. Let $G$ be a $\theta$-stable open subgroup of $(U\c)^\sigma$.
Then $G$ is a compatible real form of $U\c$ and
  $\lieu=\liek\oplus i\liep$.  
  \end{itemize} 
\end{example}

\begin{remark} \label{reductive closure remark} Let $G$ be a compatible subgroup of $U\c$.
Then the smallest complex subgroup of $U\c$ which contains $G$ is compatible and
$\theta$-stable, hence is a reductive algebraic subgroup of $U\c$.

\end{remark}
\section{Moment map decomposition}

Let $Z$ be a holomorphic $U\c$-space and $G=K\exp\liep$ a compatible
subgroup of $U\c$.  We assume that $Z$ is K\"ahler with $U$-invariant
K\"ahler structure $\omega$ and that there is a $U$-equivariant moment
mapping $\mu\colon Z\to \lieu^*$.  For any linear subspace $\liem$ of
$\lieu$ the inclusion gives by restriction a map $\mu_\liem\colon Z\to
\liem^*$. Thus we have an equivariant moment mapping $\mu_\liek\colon
Z\to \liek^*$ with respect to the $K$-action and a $K$-equivariant
mapping $Z\to (i\liep)^*$ which we denote by $\muip$.

For $\beta \in\lieu^*$ let $\m(\beta)$ denote $\mu\inv(\beta)$ and set
$\m:=\m(0)$. If $\beta\in (i\liep)^*$ then we set
$\mip(\beta):=\muip\inv(\beta)$, $\mip=\mip(0)$ and similarly for the
$\liek$ component of $\mu$.  Then $\mip$ and $\m_\liek$ are $K$-stable
and $\m$ is $U$-stable.  When necessary for clarity, we will also use
the notation $\m(Z)$, $\mip(Z)$, etc. If $\lieu=\liek+i\liep$, then
$\m=\mip\cap\m_\liek$.  Since the $U\c$-action is holomorphic, for
every $\xi\in \lieu\c$ the one-parameter group $(t,z)\mapsto(\exp
it\xi)\cdot z$, $t\in\R$, $z\in Z$, has derivative the vector field
$J\xi_Z=(i\xi)_Z$.  For a linear subspace $\lie{m}$ of $\lieu\c$ and
$z\in Z$ we set $\lie{m}\cdot z:=\{\xi_Z(z);\ \xi\in\lie{m}\}\subset
T_zZ$ and $\liem_z$ will denote $\{\xi\in\liem;\ \xi_Z(z)=0\}$.  
If $z\in Z$ is a smooth point, then, in
$T_zZ$, we use
$^{\perp_\omega}$ to denote  perpendicularity  with respect to $\omega$ and we use $^\perp$
to denote  perpendicularity with respect to the underlying Riemannian structure. If
$M\subset U\c$, then we set $M\cdot z:=\{m\cdot z;\ m\in M\}$.

\begin{lemma} \label{elementary}
Let $z\in Z$ be a smooth point. Then
\begin{enumerate}
\item $\ker d\muip(z)=(\liep\cdot z)^\perp=(i\liep\cdot z)^{\perp_\omega}$ and
\item $(\lieg\cdot z)^\perp=(\liek\cdot z)^\perp\cap\ker d\muip(z)$.
\end{enumerate}
\end{lemma}

\begin{proof} 
 This follows from the basic equation $d\mu^\xi=\imath_{\xi_Z}\omega$ and from the fact that
$\omega(J\cdot,\cdot)$ is the   Riemannian metric on $T_zZ$. 
\end{proof}

\begin{lemma} \label{zero level tangent}
Let $z\in Z$ be a smooth point and assume that $z\in\mip$. Then 
\begin{enumerate} 
\item $\lieg\cdot z=\liek\cdot z\oplus\liep\cdot z$ where 
$\liek\cdot z\perp\liep\cdot z$. 
\item If $G$ is a real form of $U\c$, then $\dim_\R\lieu\cdot z\le\dim_\R\lieg\cdot z$.
 \end{enumerate}
\end{lemma}

\begin{proof}
Since $z\in \mip$ and $\muip$ is $K$-equivariant, we have 
$\liek\cdot z\subset\ker d\muip(z)$ and therefore
$\liek\cdot z \perp \liep\cdot z$, giving (1). If $G$ is a real form of $U\c$, then
$\lieu\cdot z=\liek\cdot z +(i\liep)\cdot z$. Now $(i\liep)\cdot z=J(\liep\cdot z)$ and $\liep\cdot
z$ have the same dimension, so (2) follows from (1).
\end{proof}
 
\begin{example}  
Let $G$ be a real form of $U\c$. Assume that
$Z$ is  compact and
$U$-homogeneous, e.g., $Z$ is a flag manifold of $U\c$. 
If
$z\in \mip$, then  $\dim_\R Z=\dim_\R \lieu\cdot z\le \dim_\R\lieg\cdot z$, hence $G\cdot z$ 
is open in $Z$.

\end{example}

\begin{lemma}\label{basic equation} 
  For $z\in Z$ and $\xi\in \liep$ we have $\muip(\exp\xi\cdot
  z)=\muip(z)$ if and only if $\xi\in \liep_z$.  In particular,
  $\exp\,\liep\cdot z\cap \mip(\beta)=\{z\}$ for all $z\in
  \mip(\beta)$.
\end{lemma}

\begin{proof}
  Assume that $\muip(\exp\xi\cdot z)=\muip(z)$. Set $Y=U\c\cdot z$ and
  let $\omega_Y$ denote the K\"ahler form induced on $Y$. Let
  $\alpha(t)=(\exp t\xi)(z)$ and $\beta(t)=\mu^{i\xi}(\alpha(t))$.
  Then $\beta(0)=\beta(1)$ and
  $\beta'(t)=d\mu^{i\xi}(\xiz)(\alpha(t))=
  \omega_Y(J\xiz,\xiz)(\alpha(t))\geq 0$ since
  $\omega_Y(J\cdot,\cdot)$ is the underlying Riemannian metric of the
  K\"ahler metric on $Y$.  We must have $\beta'(0)=0$ which implies
  that $\xiz(z)=0$.  So $\xi\in\liep_z$ and $(\exp \xi)\cdot z=z$.
\end{proof}

\begin{lemma} \label{isotropy decomposition G} Let $z\in\mip$. Then
\begin{enumerate}
\item $G\cdot z\cap\mip=K\cdot z$.
\item $G_z=K_z\exp\liep_z\simeq K_z\times \liep_z$.
\end{enumerate}
In particular, $G_z$ is compatible with the Cartan decomposition of
$U\c$.
\end{lemma}

\begin{proof} Let $g=k\exp \xi\in G$ where $k\in K$ and $\xi\in\liep$, and 
  suppose that $gz\in\mip$. Then $\muip((k\exp\xi)\cdot
  z)=\muip((\exp\xi)\cdot z)=\muip(z)=0$.  Applying Lemma~\ref{basic
    equation} we obtain (1).  For (2) just notice that
  $(k\exp\xi)\cdot z=z$ implies that $\xi\in\liep_z$ by the argument
  above, so that $(k\exp\xi)\cdot z=kz=z$, and $k\in K_z$. So
  $G_z=K_z\exp \liep_z$.
\end{proof}

\begin{remark} Applying Lemma~\ref{isotropy decomposition G} 
  in the case $G=U\c$ and $z\in \m$ we obtain:
\begin{enumerate}
\item $U\c\cdot z\cap\m=U\cdot z$.
\item $(U\c)_z=U_z\exp(i\lieu_z)\simeq U_z\times i\lieu_z$.
\end{enumerate}
\end{remark}

\section{Compact isotropy groups}

Let  $G$  be a compatible Lie subgroup of $U\c$.  
Let
 $Z$ be a
holomorphic $U\c$-space with $U$-invariant K\"ahler form $\omega$
and $U$-equivariant moment mapping $\mu\colon Z\to \lieu^*$.

\begin{proposition}\label{compact isotropy}
  Let $z\in Z$ be a smooth point. If $G_z$ is compact, then $\liep_z=\{0\}$. Moreover,
  the following are equivalent.
\begin{enumerate}
\item $\liep_z=\{0\}$.
\item $d\muip(z)$ maps $\liep\cdot z$ isomorphically onto
  $(i\liep)^*$.
\item $d\muip(z)$ is surjective, i.e., $\muip$ is a submersion at $z$.
\end{enumerate}
\end{proposition}

\begin{proof}
  The Cartan decomposition of $G$ implies that $\exp\liep_z$ is
  compact if and only if $\liep_z=\{0\}$.  Thus compactness of $G_z$
  implies that $\liep_z=\{0\}$. Since $\Ker
  d\mu_{i\liep}(z)=(\liep\cdot z)^\perp$, we see that $\liep\cdot z$
  is mapped isomorphically onto $i\liep^*$ if and only if
  $\liep_z=\{0\}$ if and only if $d\muip(z)$ is surjective.
\end{proof}

\begin{remark}  Since $U\c$ is a Stein manifold,
  the isotropy group $(U\c)_z$ is compact if and only if it is finite.
  If $(U\c)_z$ is finite then $d\mu(z)$ maps $i\lieu\cdot z$
  isomorphically onto $\lieu^*$. The converse is false. Just consider
  the standard actions of $\SU(2,\C)\subset \SL(2,\C)$ on $\C^2$. Then
  the standard moment mapping has surjective differential on
  $\C^2\setminus\{0\}$, while $\SL(2,\C)$ has one dimensional isotropy
  groups.
\end{remark}

\begin{corollary} \label{corollary:compact isotropy}
  Let $G$ be closed in $U\c$ and $Z$ smooth. Then the set of $z\in \mip$ such that
  $G_z$ is compact is open in $\mip$.
\end{corollary}

\begin{proof} For $z\in\mip$ the isotropy group $G_z$ is compact
  if and only if $\liep_z=\{0\}$ (Lemma~\ref{isotropy decomposition
    G}).  Since $Z$ is smooth, the set of $z\in Z$ such that $\muip$
  has maximal rank at $z$ is open.  Now the claim follows from
  Proposition~\ref{compact isotropy}.
\end{proof}

\begin{corollary} \label{corollary compact isotropy}
Let $X$ be a
  (real) smooth $G$-stable submanifold of $Z$ and assume that $G_x$ is
  compact for all $x\in X$.  Then $\mip(\beta)\cap X$ is smooth for
  every $\beta\in (i\liep)^*$. If $x\in\mip(\beta)\cap X$, then
  $T_x(\mip(\beta)\cap X)=\Ker d(\muip |_X)_x=(\liep\cdot x)^\perp$,
  the perpendicular being taken in the tangent space $T_x(X)$.
\end{corollary}

\begin{proof} Proposition~\ref{compact isotropy} applied to 
  $U\c\cdot x$ for $x\in X$ shows that $d(\muip|_X)_x$ maps
  $\liep\cdot x$ isomorphically onto $(i\liep)^*$. Hence
  $\mip(\beta)\cap X$ is smooth and $T_x(\mip(\beta)\cap X)=\Ker
  d(\muip |_X)_x=(\liep\cdot x)^\perp$.
\end{proof}

\begin{remark} 
 Assume that $Z$ is smooth and that the $U\c$-action on $Z$ has compact (hence finite) isotropy groups. Then
  Corollary~\ref{corollary compact isotropy} says that $\m(\beta)$ is
  a smooth submanifold of $Z$ with tangent space $T_z(\m(\beta))=\Ker
  d\mu_z=(i\lieu\cdot z)^\perp$ at $z\in\m(\beta)$.
\end{remark}

Let $G$ be  
closed in $U\c$ and let $Z$ be smooth.
Let $Z_G^r$ denote $\{z\in Z;\ 
\muip \text{ is a submersion at } z\}$.  The set $Z_G^r$ is $K$-stable
and open in $Z$ (Proposition~\ref{compact isotropy}).  Let
$\Comp_{i\liep}(Z):=\{z\in Z;\ G\cdot z\cap\mip\cap
Z^r_G\ne\emptyset\}$.  We have the following slice theorem.

\begin{theorem} \label{compact slice theorem}
\begin{enumerate} 
\item The set $\Comp_{i\liep}(Z)$ is $G$-stable and for every $z\in
  \Comp_{i\liep}(Z)$ the isotropy group $G_z$ is compact. In
  particular, $\Comp_{i\liep}(Z)\subset Z^r_G$.
\item Let $S:=\mip\cap Z_G^r$. Then $S$ is a smooth closed
  $K$-submanifold of $Z_G^r$, $G\cdot S=\Comp_{i\liep}(Z)$ is open in
  $Z$ and the natural map $\twist{G}{K}{S}\to \Comp_{i\liep}(Z)$ is an
  isomorphism of $G$-manifolds.
\end{enumerate}
\end{theorem}

\begin{proof}
  The set $S=Z^r_G\cap\mip=(\muip|Z_G^r)\inv(0)$ is a closed
  submanifold of $Z_G^r$. The natural map
  $\Phi\colon\twist{G}{K}{S}\to Z$, $[g,s]\mapsto g\cdot s$, has image
  $\Comp_{i\liep}(Z)$. We first show that $\Phi$ is injective. Let
  $g_j\in G$ and $s_j\in S$ be such that $g_1\cdot s_1=g_2\cdot s_2$
  and set $g:=g_2\inv g_1$.  Then $g=k\exp\xi$ for some $k\in K$ and
  $\xi\in \liep$ and $s_2=g\cdot s_1$. This implies that $\exp\xi\cdot
  s_1=k\inv\cdot s_2\in \mip$ and therefore that $\exp\xi\cdot
  s_1=s_1$. But $G_{s_1}=K_{s_1}$ and $\liep_{s_1}=\{0\}$ by
  assumption (see Lemma~\ref{isotropy decomposition G}).
  Consequently, $\xi=0$ and $g=k$. Hence $[g_1,s_1]=[g_2 k, s_1]=[g_2,
  s_2]$. This shows injectivity of $\Phi$.
  
  In order to show that $\Phi$ is an open embedding it is sufficient
  to prove that $\Phi$ is a submersion. But for any $s\in S$ the
  tangent space $T_{[e,s]}(\twist{G}{K}{S})$ is mapped onto
  $\lieg\cdot s+(\liep\cdot s)^\perp = \liep\cdot s+(\liep\cdot
  s)^\perp=T_s(Z)$. Hence $\Phi$ is a submersion at any point $[e,s]$.
  Since $\Phi$ is $G$-equivariant it is a submersion at any point of
  $\twist{G}{K}{S}$.  Thus $\Phi$ is an isomorphism onto
  $\Comp_{i\liep}(Z)$.
\end{proof}

\begin{remark} If $S_0$ is a differentiable slice at $s_0\in S $ for 
  the $K$-action on $S$, then $K\cdot S_0$ is a slice for the
  $G$-action on $Z$ at $s_0$, i.e., $G\cdot S_0$ is open in $Z$ and
  $\twist{G}{K_{s_0}}{S_0}\to G\cdot S_0$ is a diffeomorphism.
\end{remark} 

\begin{remark} 
  The $G$-action on $\Comp_{i\liep}(Z)$ is proper.
\end{remark}

\section{Invariant plurisubharmonic functions}
\label {invariant plurisubharmonic functions}

Let $Z$ be a holomorphic $U \c$-space. Assume that we are given a
smooth $U$-invariant strictly plurisubharmonic function $\rho$ on $Z$.
Set $\mu^\xi(z)=\frac{d}{dt}|_{t=0}\rho((\exp it\xi)\cdot z)$, $z\in
Z$, $\xi\in\lieu$. Note that $\omega_Y=-dd^c(\rho|Y)$ is a
$U$-invariant K\"ahler form on every smooth complex $U$-submanifold
$Y$ of $Z$.

\begin{lemma}\label{mu lemma} Let $\rho$ and $\mu$ be as above. Then 
  $\omega=\{\rho\}$ defines a $U$-invariant K\"ahler structure on $Z$
  and $\mu$ is a moment mapping.
\end{lemma}

\begin{proof}   Let $\xi\in\lieu$ and let $Y$ be a complex 
  $U$-submanifold of $Z$. Then, since $\rho$ is $U$-invariant and
  since $U$ acts as complex analytic isomorphisms of $Z$, we have that
  $L_{\xiz}d^c(\rho|Y)=0$ where $L_{\xiz}$ denotes Lie derivative with
  respect to $\xiz$. By Cartan's formula, $L_{\xi_Z}=d\circ
  \imath_{\xi_Z}+\imath_{\xi_Z}\circ d$, hence $d\mu^\xi=d
  \imath_{\xiz}(d^c(\rho|Y))=-\imath_{\xiz}dd^c
  (\rho|Y)=\imath_{\xiz}\omega_Y$.
\end{proof}

In the situation of Lemma~\ref{mu lemma} we will say that the moment
map is \emph{defined by} or \emph{associated with} $\rho$.

\begin{remark} Let $G=K\exp\liep\subset U\c$ be compatible,  
  let $\mu$ be associated to $\rho$ and let $z\in Z$.  Since
  $\rho(K\cdot z)=\rho(z)$, $d\rho(z)(\liek\cdot z)=0$ and
  $\mip=\{z\in Z;\ \rho|{G\cdot z}\ \text{has a critical point at}\ 
  z\}$.
\end{remark}

Let $X$ be a topological space. A function $f\colon X\to \R$ is said
to be an \emph{exhaustion}, if for every $r\in\R$ the set $\{x\in X;\ 
f(x)\le r\}$ is compact.

\begin{example} Let $Z$ be a representation space of $U\c$. Let 
  $\rho(z)=||z||^2$ where $||\cdot||$ is the $U$-invariant norm on $Z$
  coming from a $U$-invariant hermitian inner product. Then $\rho$ is
  a $U$-invariant strictly plurisubharmonic exhaustion of $Z$. Every
  fiber of the quotient $\pi\colon Z\to\quot Z{U\c}$ intersects $\m$
  in a single $U$-orbit $U\cdot z$, where $U\c \cdot z$ is closed.
  Moreover, the inclusion $\m\to Z$ induces a homeomorphism
  $\m/K\simeq \quot Z{U\c}$. See \cite{KempfNess} and \cite{Schwarz}.
\end{example}

\begin{example} Let $U=\SU(2,\C)$ with its action on the complex 
  binary forms $Z$ of degree 3. Then $U\c=\SL(2,\C)$. Let
  $\rho(z)=||z||^2$ as above. The open set of orbits with finite $U\c$
  isotropy group consists of the closed $U\c$-orbits with isotropy
  group isomorphic to $\Z/3\Z$ and a non-closed orbit $U\c \cdot z_0$
  where $(U\c)_{z_0}=\{e\}$. However, $U\c\cdot z_0$ does not
  intersect $\m$, so that $\Comp_{i\liep}(Z)$ is the open set of
  closed orbits with $(\Z/3\Z)$-isotropy. This is unavoidable since
  the $U\c$-action on $\Comp_{i\liep}(Z)$ is proper. If the action
  were proper on the whole open set of orbits with finite isotropy,
  then the slice theorem for proper actions (see \cite{Palais}) would
  show that there is an open set of orbits with trivial isotropy,
  which is not the case.
\end{example}

\begin{proposition} \label{psh} Let $G\subset U\c$ be a 
  compatible subgroup and let $\rho$ be a smooth $U$-invariant
  strictly plurisubharmonic function on $Z$. Let $z\in Z$ and suppose
  that
\begin{enumerate}
\item $\rho|G\cdot z$ is an exhaustion, i.e., the map $\tilde
  \rho\colon G/G_z\to \R$, $gG_z\mapsto\rho(g\cdot z)$, is an
  exhaustion.
\end{enumerate}
Then
\begin{enumerate}
  \addtocounter{enumi}{1}
\item $\rho|G\cdot z$ has a minimum value.
\item $\rho|G\cdot z$ has a critical point.
\item $G\cdot z\cap \mip\neq\emptyset$.
\item $G\cdot z$ is closed.
\end{enumerate}
\end{proposition}
\begin{proof} Clearly we have that (1) implies (2) implies (3) 
  implies (4).  Let $(g_n\cdot z)$ be a sequence in $G\cdot z$ such
  that $\lim g_n\cdot z=z'\in Z$. Since $\{g\cdot z\in G\cdot z;\ 
  \rho(g\cdot z)\le \rho(z')+1\}$ is compact, passing to a
  subsequence, we may assume that $\lim g_n\cdot z=g_0\cdot z$ for
  some $g_0\in G$. Then $g_0\cdot z=z'\in G\cdot z$ and we have (5).
\end{proof}

\begin{remark} \label{Azad Loeb remark}
  We will show later that $(3)$ or, equivalently, $(4)$ imply that $
  \rho|G\cdot z$ is an exhaustion. Thus (1) through (4) are
  equivalent, and they all imply (5). Of course, if $\rho$ is an
  exhaustion, then (5) implies (1).  All this was first observed by
  Azad and Loeb using results of Mostow. We will obtain Mostow's
  results from general moment map properties and will then repeat the
  argument of Azad and Loeb (Lemma~\ref{spsh exhaustion}).
\end{remark} 

\begin{remark} From Lemma~\ref{isotropy decomposition G} we see 
  that the intersection in (4) is a single $K$-orbit consisting of the
  points where $\rho$ takes on its minimal value.
\end{remark}

\begin{remark} A 
smooth strictly plurisubharmonic $U$-invariant 
  exhaustion function on $Z$ exists if and only if $Z$ is a Stein
  space.
\end{remark}

\begin{example} \label{Lie algebra example} Let $\lieg$ be a real 
  semisimple Lie algebra and $G$ the corresponding adjoint group. Let
  $\lieg\c=\lieg\otimes \C$ be the complexification of $\lieg$ with
  corresponding adjoint group $G\c$.  Then $G\c$ has a Cartan
  involution $\theta$ defining a maximal compact subgroup $U$ of
  $G\c$, and $G$ is a real form of $G\c$ with compatible Cartan
  decomposition $G=K\exp\liep$.  Let $\kappa$ denote the Killing form
  on $\lieg\c$ and set $\rho(z)=-\kappa(z,\theta(z))$ for $z\in
  Z:=\lieg\c$. The function $\rho$ is a $U$-invariant strictly
  plurisubharmonic exhaustion function on $Z$. Let $\mu$ denote the
  associated moment mapping.  A simple calculation shows that
  $\m=\{z\in Z;\ [z,\theta(z)]=0\}$.  Moreover, we have $\lieg\cap
  \m=\lieg\cap\mip$ since $\lieg$ is a Lagrangian subspace of
  $\lieg\c$ with respect to $\omega=-dd^c\rho$.  In particular, as is
  well-known, the orbit $G\cdot x$ in $\lieg$ is closed if and only if
  $G\c\cdot x$ is closed in $\lieg\c$, and the latter is the case if
  and only if $x$ is semisimple.
  
  For $x \in \lieg $ we can write $x=x_{\liek}+x_\liep$ where
  $x_{\liek}\in\liek$ and $x_\liep\in\liep$.  We have $\mip\cap\lieg=
  \{x=x_\liek+x_\liep;\ [x_\liek,x_\liep]=0\}$.  Now consider the set
  $X:=\{x\in \lieg;\ G_x \text{ is compact } \}$.  Using the Jordan
  decomposition for $x\in \mip\cap X\subset \lieg\subset \lieg\c$ one
  sees that compactness of $G_x$ forces $x$ to be semisimple. Then
  $[x_\liek,x_\liep]=0$ implies that $x_\liep\in \lieg_x$, hence
  $x_{\liep}=0$ and $x\in \liek$. This shows that $X=G\cdot S$ with
  $S:=X\cap \liek$.  Applying Theorem~\ref{compact slice theorem} we
  see that $X$ is open in $\lieg$, that the natural map $\twist G K
  S\to X$ is a real analytic isomorphism and hence that $G$ acts
  properly on $X$.
\end{example}

\section{Proper actions}

As in \S 6, we consider proper $G$-actions.  
However, we do
not assume that  the  holomorphic  $U\c$-space under consideration is smooth, but we do
assume that the moment mapping is associated to a strictly plurisubharmonic
exhaustion.

Let $Z$ be a holomorphic $U\c$-space with a $U$-invariant K\"ahler
structure and moment map $\mu\colon Z\to \lieu^*$. Let $G=K\exp \liep$
be a compatible Lie subgroup of $U\c$.

\begin{proposition} \label{real proper}
  Assume that $\mu$ is associated with a strictly plurisubharmonic
  $U$-invariant exhaustion and let $X$ be a $G$-stable closed subset
  of $Z$.  If $G\times X\to X$ is a proper action, then the canonical
  map $\Phi\colon\twist{G}K{(\mip\cap X)}\to X$, $[g,x]\mapsto g\cdot
  x$, is a homeomorphism.  If $X$ is a real analytic submanifold and
  $\mu$ is real analytic, then $\Phi$ is a real analytic isomorphism.
\end{proposition}
 
\begin{proof}
  Let $x\in X$. Then $G\cdot x$ is closed, and by
  Proposition~\ref{psh} and Lemma~\ref{isotropy decomposition G},
  $G\cdot x\cap\mip=K\cdot x_0$ for some $x_0\in\mip$.  Thus $\twist
  GK{(\mip\cap X)}\to X$ is a bijection.  In order to prove that
  $\Phi$ is a homeomorphism it is sufficient to show that
  $G\times(\mip\cap X)\to X$, $(g, x)\mapsto g \cdot x$ is a proper
  map.  Let $(g_\alpha, x_\alpha)\in G\times (\mip\cap X)$ be a
  sequence such that $g_\alpha\cdot x_\alpha\to y_0\in X$.  Since
  $\rho(g\cdot x)\geq\rho(x)$ for $x\in\mip$ and $g\in G$, we have
  $\rho(x_\alpha)\le\rho(g_\alpha\cdot x_\alpha)\le r$ for some $r >
  \rho(y_0)$.  Since $\rho$ is an exhaustion, passing to a
  subsequence, we may assume that $x_\alpha\to x_0$.  Since the
  $G$-action is proper, passing to a subsequence, we may assume that
  $g_\alpha\to g_0$.  Then $y_0= g_0\cdot x_0$ and $\Phi$ is proper,
  hence a homeomorphism.
  
  If $X$ is a manifold and $\mu$ is real analytic, then $\mip\cap X$
  is smooth and $\Phi$ has maximal rank (Proposition~\ref{compact
    isotropy}).  Hence $\Phi$ is an isomorphism of manifolds.
\end{proof}

See Example \ref{Lie algebra example} for the case of $G$ acting on
the points in its Lie algebra with compact isotropy group.

\begin{corollary}\label{complex proper}
  If $G$ acts properly on $Z$, then $\twist GK\mip\to Z$ is a
  $G$-equivariant homeomorphism which is smooth for smooth $Z$.
\end{corollary}

Assume now that $G$ is a linear semisimple real algebraic group.  This
implies that $G$ is compatible with the Cartan decomposition of
$G\c=U\c=U\cdot \exp i\lieu$ where $U$ is a maximal compact subgroup
of $G\c$. Under this assumption  we can establish the real analytic
version of a theorem of Abels~\cite{Abels},
\cite{HeinznerHuckleberryKutschebauch}.
\begin{theorem} \label{proper actions}
  Let $X$ be a real analytic manifold with a proper real analytic
  $G$-action where $G$ is as above. Then there is a closed $K$-stable
  real analytic submanifold $S$ of $X$ such that the map $\twist G K
  S\to X $ is a $G$-equivariant real analytic isomorphism.
\end{theorem} 

\begin{proof}
  By \cite{HeinznerBull} there is a Stein $G\c$-manifold $Z$ and a
  closed $G$-equivariant embedding $\imath:X\to Z$. Now choose a
  $U$-invariant strictly plurisubharmonic real analytic exhaustion
  $\rho:Z\to \R$. Then we can apply Proposition~\ref{real proper}.
\end{proof} 

\begin{remark} Let $G$ be an arbitrary closed subgroup of a 
  semisimple Lie group $\hat G$ where $\hat G$ has finitely many
  components and maps injectively into its universal complexification.
  Assume that $G$ acts properly and real analytically on the real
  analytic manifold $X$. The map $x\mapsto [e,x]$ realizes $X$ as a
  $G$-stable closed subspace of $\hat X:=\twist {\hat G}GX$. Since the
  $\hat G$-action on $\hat X$ is proper we can apply
  Theorem~\ref{proper actions} to obtain a $G$-equivariant real
  analytic map $q\colon X\to \hat G/\hat K$ where $\hat K$ is a
  maximal compact subgroup of $\hat G$. Finding a global slice for the
  $G$-action on $X$ reduces to finding a global slice for the
  $G$-action on $\hat G/\hat K$.
\end{remark} 

\section{Decompositions of homogeneous spaces}
\label{decomps of homogeneous spaces}

In this section we consider proper actions on $U\c$.
We obtain Mostow
decompositions of homogeneous spaces of real reductive groups (see
\cite{Mostow1, Mostow2}).

In the following we will identify $U\c$ with $U\times i\lieu$ as a
$(U\times U)$-space (see Section \ref{Cartan decomposition}).  Let $B$
be an $(\Ad U)$-invariant inner product on $\lieu$. Define $\rho\colon
U\c\to\R$ by $\rho(u\exp(i\eta))=\frac{1}{2}B(\eta,\eta)$ for $u\in
U$, $\eta\in \lieu$.  Then $\rho$ is $(U\times U)$-invariant, and it
is a strictly plurisubharmonic function on $U\c$ (see \cite{AzadLoeb3}
or \cite{HeinznerHuckleberry}).

\begin{lemma}\label{B lemma}
  Let $\rho$ be as above. Then
\begin{enumerate}
\item $\rho$ is an exhaustion of $U\c$.
\item Let $\mu\colon U\c \to \lieu^*$ be defined using $\rho$ and the
  right action of $U$ on $U\c$. Then
  $\mu^{\xi}(u\exp(i\eta))=B(\xi,\eta)$ for $u\in U$, $\xi$, $\eta\in
  \lieu$.
\end{enumerate}
\end{lemma}

\begin{proof} 
  For (1) it is enough to show that $B$ restricted to $\lieu$ is an
  exhaustion, which is obvious. The proof of part (2) is slightly more
  complicated. Let $\xi$, $\eta\in \lieu$ and $u\in U$. Then
  $$
  \mu^\xi(u\exp(i\eta))=\frac d{dt} |
  _{t=0}\rho(u\exp(i\eta)\exp(-it\xi))=\frac d{dt} |
  _{t=0}\rho(\exp(i\eta)\exp(-it\xi))=\mu^\xi(\exp(i\eta)).
  $$
  Write $\lieu=\lieu_0\oplus\lieu_1\oplus\dots\oplus\lieu_r$ where
  each $\lieu_j$ is irreducible and the decomposition is orthogonal
  with respect to $B$.  We have $[\lieu_j,\lieu_k]=0$ for $j\ne k$.
  Write $\xi=\sum_j\xi_j$ where $\xi_j\in\lieu_j$, $j=0,\dots,r$, and
  similarly for $\eta$. Then $\mu^\xi(\exp(i\eta))=\sum_j
  \mu^{\xi_j}(\exp(i\eta))$.  Let $u_j(t)\in \exp(\lieu_j)$ and
  $\alpha_j(t)\in \lieu_j$, $t\in \R$, such that
  $\exp(i\eta_j)\exp(-it\xi_j)=u_j(t)\exp(i\alpha_j(t))$.  This gives
  $\exp(i\eta)\exp(-it\xi_j)=u_j(t)\exp(i(\alpha_j(t)+\tilde\eta_j))$
  where $\tilde\eta_j=\sum_{k\neq j} \eta_k$.  Then
  $$
  \rho(\exp(i\eta)\exp(-it\xi_j))=\frac 12
  B(\alpha_j(t),\alpha_j(t)) +\frac 12 B(\tilde\eta_j, \tilde\eta_j),
  $$
  and $\mu^\xi(\exp(i\eta))=\sum_j\mu^{\xi_j}(\exp(i\eta_j))$.
 
  From the above it is enough to prove (2) under the hypothesis that
  $\xi$, $\eta\in\lieu_0$.  Choose an embedding $ U\to U(n,\C)$ so
  that the polar decomposition of $U\c$ is induced by that of
  $U(n,\C)$.  Note that the restriction of $B$ to $\lieu_0$ is a
  positive constant times the trace form $X$, $Y\mapsto \tr(iX,iY)$,
  $X$, $Y\in\lieu\subset\lieu(n,\C)$. Thus we may assume that $B$ is
  the trace form.
  
  We have $\exp(i\eta)\exp(-it\xi)=u(t)\exp(i\alpha(t))$ where
  $u(t)\in \exp(\lieu_0)$, $\alpha(t)\in\lieu_0$ and
  $$
  \exp(2i\alpha(t))=\exp(-it\xi)\exp(2i\eta)\exp(-it\xi)=:\beta(t).
  $$
  Here $\exp$ is now the usual exponential map on matrices.  It
  suffices to establish (2) for $\xi$ and $\eta$ close to $0$ since
  the functions involved are real analytic. So we can assume that
  $2i\alpha(t)= \log\beta(t)$ is given by the usual power series
  $\log(A)=(A-I)-\frac 12(A-I)^2+\frac 13(A-I)^3\dots$ where $I$ is
  the $n\times n$ identity matrix and $A$ is near $I$. Then
  $$
  \mu^\xi(\exp(i\eta))=-\frac 12\frac d{dt}|_{t=0}\tr(\frac
  12\log\beta(t), \frac 12\log(\beta(t))=-\frac 14 \tr(\frac
  d{dt}|_{t=0} \log\beta(t), \log\beta(t))
  $$
  where $\tr(\frac d{dt}\log\beta(t),\log\beta(t))$ is a convergent
  sum of terms
  $$
  (-1)^{n+1}\frac 1n \tr[\beta'(t)(\beta(t)-I)^{n-1}+
  (\beta(t)-I)\beta'(t)(\beta(t)-I)^{n-2}+
  \dots+(\beta(t)-I)^{n-1}\beta'(t),\log\beta(t)].
  $$
  Using the identities satisfied by trace and the fact that
  $\beta(t)$ and $\log\beta(t)$ commute we can permute the terms
  involving $\beta'(t)$ and $\beta(t)-I$ to obtain
  $$
  -\frac 14 \tr(\frac d{dt}\log\beta(t), \log\beta(t))=-\frac 14
  \tr(\beta'(t)\beta(t)^{-1},\log\beta(t)).
  $$
  Evaluating at $t=0$ and using the identities satisfied by trace
  we obtain
  $$
  \mu^\xi(\exp(i\eta))=-\frac
  14\tr([-i\xi\exp(2i\eta)+\exp(2i\eta)(-i\xi)]\exp(-2i\eta),
  2i\eta)=\tr(i\xi,i\eta)=B(\xi,\eta).
  $$
\end{proof}

\begin{remark}
  If we identify $\lieu$ and $\lieu^*$ using the inner product $B$, we
  get that $\mu(u\exp(i\eta))=\eta$ for all $u\in U, \eta\in \lieu$.
  Thus after identifying $U\c$ with $U\times i\lieu\cong U\times
  \lieu$ the moment map is given by projection on the second
  component! Also, if we view $U\c\cong U\times \lieu^*$ as the
  cotangent bundle of $U$, then $\mu$ coincides with the standard
  moment map on the cotangent bundle associated with the standard
  symplectic form. Moreover, the complex structure induced on the
  tangent bundle of $U$ by $U\c\cong U\times \lie u$ is the so called
  adapted complex structure with respect to the Riemannian metric on
  $U$ defined by $B$ (see \cite{LempSz} and \cite{GuSt}).
\end{remark}  

Let $G=K\exp\liep$ and $H=L\exp \lieq$ be compatible closed subgroups
of $U\c$. Assume that $H$ is a subgroup of $G$ such that $L\subset K$
and $\lieq\subset\liep$.  We apply the results of the previous section
to the free proper actions of $U\c$ and its compatible subgroups on
$Z:=U\c$ by multiplication on the right; $g\cdot z=zg^{-1}$.

\begin{theorem} \label{Mostow1}  
  Write $\liep=\lieq+\lieq^\perp$ where $\lieq^\perp$ is the
  perpendicular of $\lieq$ in $\liep$ relative to $B$. Then the
  $(H\times K)$-equivariant map $H\times (K\times\lieq^\perp)\to G$,
  $(h,(k,\xi))\mapsto k\exp(\xi)h^{-1}$, induces an $(H\times
  K)$-equivariant isomorphism
  $$\twist HL{(K\times\lieq^\perp})\tosim G$$
  where
  $l(h,k,\xi)=(hl\inv,k l\inv,\Ad l(\xi))$, $l\in L$, $h\in H$, $k\in
  K$, $\xi\in q^\perp$.  In particular, we have an induced
  $K$-equivariant isomorphism
  $$\twist KL{\lieq^\perp}\to G/H\ .$$
\end{theorem}
\begin{proof} This is an application of Proposition~\ref{real proper} 
  where $X$ is $G$, $Z$ is $U\c$, $H$ plays the role of $G$ and $\rho$
  is as in Lemma~\ref{B lemma}. From Lemma~\ref{B lemma} we get the
  equivariant identification of 
$\m_{i\lieq}$
with $K\times
  \lieq^\perp$.
\end{proof}

\begin{remark} Theorem~\ref{Mostow1} implies that, as sets,
  $G=K\exp(\lieq^\perp)H$. Applying the same reasoning to $U\c$ and
  $G$ we obtain that $U^\C=U\exp(i\liek)G$.
 \end{remark}
 
 Let $M$ be a compact Lie subgroup of $U$ and $\liem$ its Lie algebra.

\begin{corollary} \label{Mostow2}  
  Define $\liem^\perp$ to be the perpendicular of $\liem$ in $\lieu$
  relative to $B$ and let $M$ act on $U\times i\liem^\perp$ by
  $g\cdot(u,i\xi)=(ug^{-1},i\mathrm{Ad}(g)\xi)$, $g\in M$, $u\in U$,
  $\xi\in\liem^\perp$.  Then the map $M\c\times (U\times
  i\liem^\perp)\to U\c$, $(h,u,i\xi)\mapsto u\exp(i\xi)h^{-1}$,
  induces an $(M\c\times U)$-equivariant isomorphism
  $$\twist {M\c} M{(U\times i\liem^\perp})\tosim U\c$$
  where $M\c$
  acts on $U\c$ by right multiplication and $M$ acts on $M\c$ by right
  multiplication. In particular, we have an induced $U$-equivariant
  isomorphism
  $$\twist UM{i\liem^\perp}\tosim U\c/M\c\ .$$
 \end{corollary}
 
 Let $Z$ be a holomorphic $U^\C$-space and $\rho\colon Z\to \R$ a
 $U$-invariant smooth strictly plurisubharmonic function with
 associated moment map $\mu\colon Z\to \lieu^*$. Let $G=K\exp\liep$ be
 a closed compatible subgroup of $U\c$. The following was proved by
 Azad and Loeb (\cite{AzadLoeb2}) using Mostow's results. For
 completeness, we reproduce their argument.

 \begin{lemma} \label{spsh exhaustion}
   The restriction $\rho|{G\cdot z}$ is an exhaustion if (and only if)
   $G\cdot z\cap \mip\ne\emptyset$.
 \end{lemma}

 \begin{proof}
   Let $z_0\in G\cdot z\cap \mip$. Then
   $G_{z_0}=K_{z_0}\exp(\liep_{z_0})$ and by Theorem~\ref{Mostow1},
   $(\liep_{z_0})^\perp\to Z$, $\xi\to \exp(\xi)\cdot z_0$, is an
   injective immersion.  Let $\tilde\rho\colon (\liep_{z_0})^\perp\to
   \R$ be defined by $\tilde\rho(\xi)=\rho(\exp(\xi)\cdot z_0)$. Since
   $\C\to \R$, $w\mapsto\rho(\exp(w\xi)\cdot z_0)$ is strictly
   subharmonic for $\xi\ne 0$ and does not depend on the imaginary
   part of $w$ (since $\rho$ is $U$-invariant), it follows that
   $\tilde \rho|{\R\xi}$ is strictly convex for every $\xi\ne 0$.
   Hence $0\in(\liep_{z_0})^\perp$ is the unique critical point of
   $\tilde\rho$. This all implies that $\tilde \rho$ is an exhaustion
   \cite{AzadLoeb1}.  It follows that $\rho|{G\cdot z}$ is an
   exhaustion.
\end{proof}

\section{Quotients by complex groups}\label{quotients by reductive groups}

Before we continue our study of $G$-actions, we recall some results
about actions of complex reductive groups.  Let $Z$ be a holomorphic
$U\c$-space with $U$-invariant K\"ahler structure $\omega$ and moment
map $\mu$.  Set $\mathcal{S}_{U\c}(\m)=\{z\in Z;\ \overline{U\c\cdot
  z}\cap \m\neq\emptyset\}$, the set of \emph{semistable points of $Z$ with
respect to $\mu$ and the $U\c$-action}.
 The
following result can be found in \cite{HeinznerLoose}.

\begin{theorem} \label{Uc quotients}
\begin{enumerate}
\item The set $\mathcal{S}_{U\c}(\m)$ is open and $U\c$-invariant, and
  there is an analytic Hilbert quotient $\pi:\mathcal{S}_{U\c}(\m)\to
  \quot{\mathcal{S}_{U\c}(\m)}{U\c}$.
\item Each fiber of $\pi$ contains a unique closed $U\c$-orbit which
  is the unique orbit of minimal dimension in that fiber.
\item For every $z\in \m$, the orbit $U\c\cdot z$ is closed in $\mathcal{S}_{U\c}(\m)$.
\item The inclusion $\m\to \mathcal{S}_{U\c}(\m)$ induces a
  homeomorphism
  $$\m/U\cong \quot{\mathcal{S}_{U\c}(\m)}{U\c}\ .$$
\end{enumerate}
\end{theorem}

\begin{remark} \label{global rho} Suppose that $\mu$ and $\omega$ are 
  associated to a smooth strictly plurisubharmonic exhaustion $\rho$
  and $z\in Z$. Then $\rho| \overline{U\c\cdot z}$ takes on a minimum
  value, hence $z\in \mathcal{S}_{U\c}(\m)$. Thus
  $Z=\mathcal{S}_{U\c}(\m)$.
\end{remark}

Regarding the existence of strictly plurisubharmonic functions we have
the following result from \cite{HeinznerLoose},
\cite{HeinznerHuckleberryLoose} and \cite{HeinznerHuckleberryInv}.

\begin{theorem} \label{local plurisubharmonic} Let $Z$, $\omega$ and $\mu$ 
  be as above, and let $z_0\in\m$. Then there is a neighborhood $Q_0$
  of $\pi(z_0)\in\quot{\mathcal{S}_{U\c}(\m)}{U\c}$ such that, on
  $\Omega_0:=\pi\inv(Q_0)$, $\mu$ and $\omega$ are associated to a
  strictly plurisubharmonic function $\rho$ such that
\begin{enumerate}
\item $(\pi\times\rho)\colon \Omega_0\to Q_0\times\R$ is proper and
  $\rho$ is bounded from below.
\item $(\pi\times\rho)\inv(\pi(z_0),\rho(z_0))=U\cdot z_0$.
\end{enumerate}
\end{theorem}

\begin{remark} Property (1) is not stated in 
  \cite{HeinznerHuckleberryInv}.  But an easy modification of the
  proof of the exhaustion theorem in \cite{HeinznerHuckleberryInv}
  gives the claimed statement (see \cite{HeinznerHuckleberry}).
\end{remark}

\begin{remark}\label{Hilbert quotient remark} Suppose that the analytic Hilbert quotient 
  $\pi\colon Z\to\quot Z{U\c}$ exists. Then there is an open cover
  $\{\Omega_\alpha\}$ of $\quot Z{U\c}$ where each $\Omega_\alpha$ and
  $\pi\inv(\Omega_\alpha)$ are Stein.  Moreover, there are strictly
  plurisubharmonic $U$-invariant exhaustions $\rho_\alpha$ on
  $\pi\inv(\Omega_\alpha)$.  Thus, locally, we are in the case that
  $Z=\mathcal{S}_{U\c}(\m)$.
\end{remark}

Finally, from \cite{HeinznerLoose} we have the existence of slices.
\begin{theorem} \label{holomorphic slice theorem} Let $z\in\m$. Then there is
  a locally closed $(U\c)_z$-invariant subspace $S\subset Z$, $z\in
  S$, such that $V:=U\c\cdot S$ is open in $Z$ and such that the
  natural $U\c$-invariant holomorphic map $\varphi\colon
  \twist{U\c}{(U\c)_z}S\to V$, $[g,s]\mapsto g\cdot s$, is
  biholomorphic.
\end{theorem}

\begin{remark} \label{smooth Stein}
  By construction of $S$, we may assume that it admits a
  $(U\c)_z$-equivariant closed embedding into an open $(U\c)_z$-stable
  Stein submanifold $\hat S$ of the $(U\c)_z$-representation space
  $T_zZ$. This implies, in particular, that $\twist{U\c}{(U\c)_z}S$
  can be realized as a closed subspace of the Stein $U\c$-manifold
  $\twist{U\c}{(U\c)_z}{\hat S}$.
\end{remark}

\section{Semistable points with respect to compatible subgroups}
\label{semistable}

Let $Z$ be a holomorphic $U\c$-space with equivariant moment map
$\mu\colon Z\to \lieu^*$ and let $G=K\exp\liep$ be a compatible closed
subgroup of $U\c$.  For a subset $Y$ of $Z$ and a subset $H$ of $U\c$
we define $\mathcal{S}_H(Y):=\{z\in Z;\ \overline{H\cdot z}\cap
Y\ne\emptyset\}$.  We call $\mathcal{S}_G(\mip)$  
the set of \emph{semistable points of $Z$ with respect
to $\muip$ and the $G$-action}.
We already met $\mathcal{S}_{U\c}(\m)$ in the
previous section.  When necessary for clarity, we will use notation
such as $\mathcal{S}_{U\c}(\m(Z))$ and $\mathcal{S}_G(\mip(Z))$.

\begin{remark} If $ \beta\in\lieu^*$ is a fixed point with respect to 
  the co-adjoint $U$-action, then $z\mapsto \mu(z)-\beta$ defines a
  shifted $U$-equivariant moment map $\mu-\beta:Z\to \lieu^*$. The set
  of semistable points with respect to this shifted moment map is
  $\mathcal{S}_{U\c}(\m(\beta))$. If $\beta$ annihilates $i\liep$,
  then $(\mu-\beta)_{i\liep}=\mu_{i\liep}$ and the set of semistable
  points with respect to $G$ remains unchanged.
\end{remark}

\begin{proposition} \label{G-orbit closure} 
  Let $G$ be a closed compatible subgroup of $U\c$.  If
  $Z=\mathcal{S}_{U\c}(\m)$, then
    \begin{enumerate}
    \item every $G$-orbit in $Z$ contains a closed $G$-orbit in its
      closure,
    \item $Z=\mathcal{S}_G(\mip)$ and
    \item the closed $G$-orbits are precisely the orbits intersecting
      $\mip$.
    \end{enumerate}
\end{proposition}

\begin{proof} 
  We may assume that $Z$ coincides with a fiber of $\pi\colon Z\to
  \quot{Z}{U\c}$.  The moment map $\mu$ is then associated to a
  strictly plurisubharmonic exhaustion $\rho$.  Now let $z\in Z$.
  Since $\rho | \overline{G\cdot z}$ is an exhaustion, it has a
  minimum at a point $z_0$. Thus $\mu_{i\liep}(z_0)=0$, so
  $z\in\mathcal{S}_G(\mip)$ and we have (2). For any $z_0\in\mip$,
  Lemma \ref{spsh exhaustion} shows that $\rho | G\cdot z_0$ is an
  exhaustion, so $G\cdot z_0$ is closed, and we have (1) and (3).
\end{proof}

\begin{remark} \label{real proper remark} In many instances, one can use Theorem~\ref{local
    plurisubharmonic} to replace the hypothesis of the existence of a
  strictly plurisubharmonic $U$-invariant exhaustion by the hypothesis
  that $Z=\mathcal{S}_{U\c}(\m)$. For example, this can be done in the
  case of Proposition~\ref{real proper} and Corollary~\ref{complex
    proper}
\end{remark} 

Now we have an important definition.

\begin{definition} \label{muip-adapted}
  Let $A\subset Z$ be $K$-stable. For $z\in Z$ and $\xi \in \lieu$ let
  $I_A(z;\xi):=\{t\in\R;\ (\exp it\xi)\cdot z\in A\}$. Then $A$ is
  said to be \emph{$\muip$-adapted\/} if for every $z\in Z$, $\xi\in
  i\liep$ and every nonempty connected component $C$ of $I_A(z;\xi)$
  the following holds.

\begin{itemize}
\item[i)] If $t_-=\inf C>-\infty$, then $\muip^{\xi}((\exp
  it_-\xi)\cdot z)<0$ and
\item[ii)] if $t_+=\sup C<+\infty$, then $\muip^{\xi}((\exp
  it_+\xi)\cdot z)>0$.
\end{itemize}
\end{definition}
\begin{remark} 
  Conditions i) and ii) are valid for every $z\in Z$ if and only if
  they are valid for every $z\in A$.
\end{remark}

\begin{remark}\label{adapted}
  Let $A$ be $\mu_{i\liep}$-adapted.  Let $\xi\in i\liep$ and $z\in
  Z$. If $\xi_Z(z)=0$, then $I_A(z;\xi)=\emptyset$ or $\R$. Assume
  that $\xi_Z(z)\neq 0$ and that $I_A(z;\xi)\neq\emptyset$. Clearly
  $I_A(z;\xi)$ is not a point. As we saw in the proof of Lemma
  \ref{isotropy decomposition G}, the function $\lambda\colon\R\to
  \R$, $\lambda(t)=\mu^{\xi}(\exp it\xi\cdot z)$, is strictly
  increasing. It follows that $I_A(z;\xi)$ cannot contain two disjoint
  components, i.e., $I_A(z;\xi)$ is connected.
\end{remark}
\begin{remark} If $G=U\c$, then $\muip=\mu$ and we have 
  the notion of $\mu$-adapted sets.
\end{remark}
We leave the proof of the following to the reader.
\begin{lemma} \label{mip union lemma} 
  Arbitrary unions and finite intersections of $\muip$-adapted sets
  are $\muip$-adapted.
\end{lemma}

Here is an important tool for constructing $\mu$-adapted neighborhoods
(c.f. \cite{Heinzner}).

\begin{proposition}  \label{adapted sets}
  Let $Z$ be a holomorphic $U\c$-space.  Assume that $\mu\colon
  Z\to\lieu^*$ is associated with a $U$-invariant strictly
  plurisubharmonic smooth function $\rho\colon Z\to \R$.  For $r\in
  \R$ set $\Delta_r(\rho):=\{z\in Z;\ \rho(z)<r\}$. Then
  $\Delta_r(\rho)$ is $\mu$-adapted.
\end{proposition}

\begin{proof}
  Let $\xi\in \lieu$ and $z_0\in \Delta_r(\rho)$. By definition we
  have
  $$
  \mu^\xi(\exp(it\xi)\cdot z_0)=\frac{d}{dt}\,
  \rho(\exp(it\xi)\cdot z_0) \,.$$
  
If $\xi_Z(z_0)=0$ we have nothing to
  prove, so assume that $\xi_Z(z_0)\ne 0$. Then the function $t\to
  \mu^\xi(\exp(it\xi)\cdot z_0)$ is strictly increasing and the
  boundary of $\Delta_r(\rho)$ is contained in $\{z\in Z;\ 
  \rho(z)=r\}$.  Since $t\mapsto\rho(\exp(it\xi)\cdot z_0)$ is a
  strictly convex function this implies that $\Delta_r(\rho)$ is
  $\mu$-adapted.
\end{proof}

The following lemma shows the utility of the notion of
$\muip$-adaptedness.

\begin{lemma} \label{intersection}
  If $A_1, A_2\subset Z$ are $\muip$-adapted, then $G\cdot (A_1\cap
  A_2)\supset A_1\cap G\cdot A_2$. In particular,
  $$G\cdot (A_1\cap A_2)=G\cdot A_1\cap G\cdot A_2\ .$$
\end{lemma}

\begin{proof}
  Suppose that $x_1 = g\cdot x_2$ where $g\in G$ and $x_j\in A_j$,
  $j=1$, $2$.  There is a $k\in K$ and $\xi\in i\liep$ such that
  $g=k\exp i\xi$. Replacing $x_1$ by $k\inv x_1$ we may assume that
  $g=\exp(i\xi)$. Consider the path $\alpha(t)=\exp(i t\xi)\cdot x_2$
  from $[0,1]$ to $Z$. It is enough to show that $\alpha(t_0)\in
  A_1\cap A_2$ for some $t_0\in[0,1]$. This is obvious if $x_1\in A_2$
  or if $x_2\in A_1$, so we may assume that there are $t_1$, $t_2\in
  [0,1]$ where $\alpha$ enters $A_1$, resp.\ leaves $A_2$. From
  $\muip$-adaptedness we get that $\mu^\xi(\alpha(t_1))<0$ and
  $\mu^\xi(\alpha(t_2))>0$. But $\mu^{\xi}$ is increasing on the image
  of $\alpha$, so $t_1<t_2$ and $\alpha(t)\in A_1\cap A_2$ for
  $t\in[t_1,t_2]$.
\end{proof}

\begin{remark} \label{locally g stable}
  Let $\Omega$ be open and $\muip$-adapted and let $A$ be a $K$-stable
  closed subset of $\Omega$. Assume that $A$ is locally $G$-stable,
  i.e, for all $a\in A$ and 
$\xi\in i\lie p$ the set $I_A(a;\xi)$ is a
  neighborhood of $0\in \R$. Since $A$ is closed in $\Omega$,
  $I_A(a;\xi)$ is closed in $I_\Omega(a;\xi)$. Now the condition that
  $A$ is locally $G$-stable implies that $I_A(a;\xi)$ is open. Thus
  $I_A(a;\xi) =I_\Omega(a;\xi)$ and consequently $A$ is
  $\muip$-adapted. The equality $I_A(a;\xi) =I_\Omega(a;\xi)$ also
  implies that $A=G\cdot A\cap \Omega$.  In particular, $G\cdot A$ is
  closed in $G\cdot \Omega$ and if $A$ is real analytic in $\Omega$,
  then $G\cdot A$ is real analytic in $G\cdot \Omega$. Note that if
  $\Omega$ is $G$-stable then local stability of $A$ implies that it
  is $G$-stable.
\end{remark}

\begin{corollary}\label{omega cor} 
  Suppose that every $K$-orbit in $\mip$ has a neighborhood basis of
  open $\muip$-adapted subsets. Let $C_1$, $C_2\subset\mip$ be closed,
  disjoint and $K$-stable. Then there are $G$-invariant open
  neighborhoods $\Omega_i$ of $C_i$, $i=1$, $2$, such that
  $\Omega_1\cap\Omega_2=\emptyset$.
\end{corollary}
\begin{proof} Let $W_1$ and $W_2$ be disjoint open 
  $K$-invariant neighborhoods of $C_1$ and $C_2$, respectively. Then
  each point $z\in C_1$ has an open $\muip$-adapted neighborhood
  $W_z\subset W_1$. Let $W_1'$ be the union of the $W_z$ for $z\in
  C_1$, and construct $W_2'$ similarly. Then $W_1'$ and $W_2'$ are
  disjoint and $\muip$-adapted, hence the sets $\Omega_i:=G\cdot
  W_i'$, $i=1$, $2$, are $G$-invariant disjoint open neighborhoods of
  $C_1$ and $C_2$.
\end{proof}

It would be ideal if we could show that every $K$-orbit in $\mip$ has
a basis of open $\mu_{i\liep}$-adapted neighborhoods. We are able to
do this in the case that $Z= \mathcal{S}_{U\c}(\m)$ (Theorem \ref{main
  theorem}).  However, we can get many results with the weaker
property established in Corollary \ref{omega cor}.  This will follow
from the special case where $U$ is commutative, which we handle in the
next section.

\begin{definition}\label{defn separation property} 
  We say that the {\it separation property holds\/} if any two
  $K$-stable disjoint closed subsets of $\mip$ are contained in
  disjoint open $G$-stable sets.

\end{definition}

\begin{proposition}\label{intersection with mip} 
  Suppose that the separation property holds and let $z\in Z$. Then
  $\overline{G\cdot z}\cap \mip\ne\emptyset$ if and only if
  $\overline{G\cdot z}\cap \mip=K\cdot z_0$ for some
  $z_0\in\overline{G\cdot z}\cap \mip$.
\end{proposition}

\begin{proof}
  Assume that $K\cdot x_1, K\cdot x_2\subset \overline{G\cdot
    z}\cap\mip$ and that $K\cdot x_1\ne K\cdot x_2$. Let $\Omega_j$ be
  $G$-invariant disjoint neighborhoods of $K\cdot x_j$, $j=1$, $2$.
  Then $G\cdot z\subset \Omega_1\cap \Omega_2=\emptyset$, a
  contradiction.
\end{proof}

\begin{definition} Let $X$ be a topological $G$-space. We say that the 
  \emph{topological Hilbert quotient of $X$ by $G$\/} exists if $\sim$
  is an equivalence relation, where $x\sim y$ for $x$, $y\in X$ if and
  only if $\overline{G\cdot x}\cap \overline{G\cdot y}\ne \emptyset$.
  If this is the case, then we define the Hilbert quotient to be the
  set of equivalence classes, denoted $\quot{X}{G}$, with the quotient
  topology. If $Y$ is another topological $G$-space with topological
  Hilbert quotient $\quot YG$ and if $\varphi\colon X\to Y$ is
  continuous and $G$-equivariant, then $\quot\varphi G$ will denote
  the induced continuous map from $\quot XG$ to $\quot YG$.
\end{definition}

\begin{remark} Let $Z$ be a holomorphic $U\c$-space. If the analytic 
  Hilbert quotient $\pi\colon Z\to\quot Z{U\c}$ exists, then it is
  also the topological Hilbert quotient since $\pi(x)=\pi(y)$ if and
  only if $\overline{U\c\cdot x}$ and $\overline{U\c\cdot y}$ contain
  the same closed $U\c$-orbit.
\end{remark}

\begin{corollary}\label{intersection with mip corollary} Suppose that the
  separation property holds. Then the topological Hilbert quotient
  $\quot{\mathcal{S}_G(\mip)}G$ exists.  We have $x\sim y$ in
  $\mathcal{S}_G(\mip)$ if and only if $K\cdot x_0= K\cdot y_0$ where
  $K\cdot x_0=\overline{G\cdot x}\cap \mip$ and $K\cdot
  y_0=\overline{G\cdot y}\cap \mip$. The inclusion of $\mip$ into
  $\mathcal{S}_G(\mip)$ induces a continuous bijection of $\mip/K$
  with $\quot{\mathcal{S}_G(\mip)}G$.
\end{corollary}

\begin{proof} Suppose that 
  $z\in\overline{G\cdot x}\cap\overline{G\cdot y}$,
  $z\in\mathcal{S}_G(\mip)$, and that $\overline{G\cdot x}\cap
  \mip=K\cdot x_0$ and that $y_0$ and $z_0$ are defined similarly.
  Then by Proposition \ref{intersection with mip} we must have that
  $K\cdot x_0=K\cdot z_0=K\cdot y_0$.
\end{proof}

\begin{remark} 
  Assume  the separation property.
  Then every fiber of the quotient map $\pi\colon \mathcal{S}_G(\mip)\to
  \quot{\mathcal{S}_G(\mip)}{G}$ intersects $\mip$ in a unique
  $K$-orbit $K\cdot z_0$.  The orbit   $G\cdot z_0$ is
closed in $\mathcal{S}_G(\mip)$ and is the unique orbit of lowest dimension in the
corresponding fiber of $\pi$ (Corollaries
\ref{closed orbits of G} and \ref{minimal dimension orbit}). Conversely, any closed orbit
in $\mathcal{S}_G(\mip)$ intersects $\mip$, by definition. Thus the quotient map
$\pi$ gives a
  parametrization of the   closed $G$-orbits in $\mathcal{S}_G(\mip)$.
\end{remark}
 
\begin{lemma} Assume the separation property. Then 
  $\pi | \mip\colon \mip\to \quot{\mathcal{S}_G(\mip)}{G}$ is a closed
  mapping.
\end{lemma}

\begin{proof} Let $C\subset\mip$ be closed. We may assume that 
  $C$ is $K$-invariant.  Set $Y:=\pi\inv(\pi(C))$. Then we must show
  that $Y$ is closed in $\mathcal{S}_G(\mip)$, so let $y_n\in Y$ and
  suppose that $y_n\to y\in\mathcal{S}_G(\mip)$. Let
  $y_0\in\overline{G\cdot y}\cap\mip$. If $y_0\not\in C$, then we can
  find disjoint open $G$-stable neighborhoods $\Omega_1$ of $y_0$ and
  $\Omega_2$ of $C$. For large $n$ we must have that $y_n\in \Omega_1$
  which implies that $\overline{G\cdot y_n}$ lies in the complement of
  $\Omega_2$. But then we cannot have $y_n\in\pi\inv(\pi(C))$, a
  contradiction. Hence $y_0\in C$, $Y$ is closed in
  $\mathcal{S}_G(\mip)$ and $\pi(C)$ is closed.
\end{proof}

\begin{corollary}  Assume the separation property. If $A$ is 
  a $G$-stable closed subset $\mathcal{S}_G(\mip)$, then $\pi(A)$ is
  closed. If $A_1$ and $A_2$ are two closed $G$-stable subsets of
  $\mathcal{S}_G(\mip)$, then we have $\pi(A_1)\cap
  \pi(A_2)=\pi(A_1\cap A_2)$.
\end{corollary}

\begin{proof} We have that $\pi(A)=\pi(A\cap\mip)$ 
  is closed. Similarly, $\pi(A_1)\cap \pi(A_2)=
  \pi(A_1\cap\mip)\cap\pi(A_2\cap\mip)=\pi(A_1\cap A_2)$.
\end{proof}

\begin{corollary} \label{complete} 
  Assume the separation property. Then the bijection
  $\mip/K\to\quot{\mathcal{S}_G(\mip)}G$ is a homeomorphism. In
  particular, $\quot{\mathcal{S}_G(\mip)}G$ is  metrizable and locally compact.
\end{corollary}

\begin{proof} 
  By the lemma, the bijection is a closed mapping.
\end{proof}

\begin{remark} \label{open map remark} 
  Note that the quotient map $\pi$ is open at every point of $\mip$
  and therefore open at every point of $G\cdot \mip$.
\end{remark} 

\section{Commutative groups}\label{commutative groups}

In this section $U$ will denote a commutative compact Lie group.  Let
$Z$ be a holomorphic $U\c$-space with moment map $\mu\colon Z\to
\lieu^*$ and let $G=K\exp\liep$ be a compatible Lie subgroup of $U\c$.
Let $\Ann i\liep$ denote the annihilator of $i\liep$ in $\lieu^*$.
 
\begin{proposition}   
  $$\mathcal{S}_G(\mip)= \bigcup_{\beta\in\Ann i\liep}
  \mathcal{S}_{U\c}(\m(\beta))\, .$$
\end{proposition}  

\begin{proof} 
  Let $z\in \mathcal{S}_{U\c}(\m(\beta))$ where $\beta\in\Ann i\liep$.
  On the fiber $Y$ through $z$ of the quotient
  $\pi_\beta\colon\mathcal{S}_{U\c}(\m(\beta)) \to
  \quot{\mathcal{S}_{U\c}(\m(\beta))}{U\c}$ the shifted moment map
  $\mu-\beta$ is associated with a strictly plurisubharmonic
  $U$-invariant smooth exhaustion $\rho\colon Y\to \R$. Since
  $\rho|\overline{G\cdot z}$ attains a minimum value we obtain that
  $\overline{G\cdot z}\cap \mip\ne\emptyset$. Consequently,
  $\mathcal{S}_{U\c}(\m(\beta))\subset \mathcal{S}_G(\mip)$.
  
  If $z\in\mathcal{S}_G(\mip)$, choose $z_0\in\mip\cap\overline{G\cdot
    z}$.  Then $z_0\in\m(\beta_0)$ where $\beta_0:=\mu(z_0)\in\Ann
  i\liep$. This gives $z_0\in\m(\beta_0)\cap\overline{G\cdot z}\subset
  \m(\beta_0)\cap\overline{U\c\cdot z}$, i.e., $z\in\mathcal
  {S}_{U\c}(\m(\beta_0))$ .
\end{proof}

\begin{remark} If $\lieu=i\liep$, then 
  $\mathcal{S}_G(\mip)=\mathcal{S}_{U\c}(\m)$.
\end{remark}
 
\begin{corollary} We have that $\mathcal{S}_G(\mip)$ is an 
  open $U\c$-stable subset of $Z$.
\end{corollary}

\begin{remark} 
  We conjecture that, even for noncommutative groups,
  $\mathcal{S}_G(\mip)$ is always open in $Z$.
\end{remark}

In the following we will assume that $G=K\exp\liep$ is a closed
compatible subgroup of $U\c$.

\begin{proposition} \label{mip-adapted sets commutative case}
  Every $K$-orbit $K\cdot z_0\subset \mip$ has a neighborhood basis
  consisting of $\muip$-adapted open sets.
\end{proposition}

\begin{proof}
  Shifting by an element in $\Ann i\liep$ we may assume that $z_0\in
  \m$. Let $M$ denote $U_{z_0}$. Then   by the holomorphic slice theorem
(Theorem
  \ref{holomorphic slice theorem}) we may assume that
  $Z=\twist{U\c}{M\c}S$ for some $M\c$-space $S$, and we
  may assume that the moment map is associated to a strictly
  plurisubharmonic function $\rho$.  Let $\pi\colon Z\to
  \quot{Z}{U\c}$ denote the quotient map and let $q\colon Z\to
  U\c/M\c\cong U\c\cdot z_0$ denote the bundle projection. We may
  assume that $\rho\times\pi$ is proper, that $\rho$ is bounded from
  below and that $(\rho\times\pi)\inv(\rho(z_0),\pi(z_0))=U\cdot z_0$.
  Moreover, we have that $U\cdot z_0\cap G\cdot z_0=K\cdot z_0$. For
  suppose that $u\cdot z_0=k\exp\xi\cdot z_0$ where $u\in U$, $k\in K$
  and $\xi\in\liep$.  Then $u\inv k\exp\xi\in (U\c)_{z_0}=U_{z_0}\exp
  i\lieu_{z_0}$, so $u\inv k\in U_{z_0}$ and $u\cdot z_0=u(u\inv
  k)\cdot z_0=k\cdot z_0$.
  
  Since $G$ is closed in $U\c$ and $z_0\in\m$, the group $G/G_{z_0}$
  acts properly and freely on $U\c\cdot z_0$.  Thus the quotient
  $U\c\cdot z_0/G$ is a Hausdorff space.  Let $\bar q\colon Z\to
  (U\c\cdot z_0)/G$ denote $q$ composed with the quotient map
  $U\c\cdot z_0\to(U\c\cdot z_0)/G$. The map $(\pi\times \rho\times
  \bar q)\colon Z\to \quot{Z}{U\c}\times \R\times (U\c\cdot z_0)/G$ is
  proper and satisfies $(\pi\times \rho\times \bar
  q)^{-1}(\pi(z_0),\rho(z_0),\bar q(z_0)) =U\cdot z_0\cap G\cdot
  z_0=K\cdot z_0$. This implies that the sets $\pi\inv(Q)\cap
  \rho^{-1}((-\infty, r))\cap \bar q^{-1}(B)$ where $Q$ is an open
  neighborhood of $\pi(z_0)$, $r>\rho(z_0)$ and $B$ is an open
  neighborhood of $\bar q(z_0)$ form a basis of open $K$-invariant
  neighborhoods of $K\cdot z_0$ in $Z$.  Moreover, $\pi\inv(Q)\cap
  \rho^{-1}((-\infty, r))$ is $\mu$-adapted, hence $\muip$-adapted,
  and $\bar q^{-1}(B)$ is $\muip$-adapted since it is $G$-invariant.
  Therefore the intersection remains $\muip$-adapted.
\end{proof}

\begin{remark} This paper makes essential use of the holomorphic slice
theorem and the existence of a potential  $\rho$ only in the case that $U$ is commutative.
One can give short proofs of the required results in this case, thus circumventing the use of
the results in \cite{HeinznerLoose}. Then our results in this paper generalize those of
\cite{HeinznerLoose} and are independent of the results there. 

\end{remark}
\section{The Separation Property} 

In this section we establish the separation property. Let $U$ be a
compact Lie group.

\begin{proposition} \label{KAK result} Let 
  $G\subset U\c$ be closed and compatible. Let $A$ be a maximal
  connected subgroup of $\exp\liep$ and let $\liea$ denote its Lie
  algebra. Then $\liea$ is Abelian, so that $A\simeq \R^l$ for some
  $l$. We have $\Ad(K)\liea=\liep$ and $G=KAK$.
\end{proposition}

\begin{proof}
  Since $\left[\liep,\liep\right]\subset \liek$, both $\liea$ and $A$
  are commutative. 
  Write $\lieg=\lieg_{s}\oplus\lieg_r$ where
  $\lieg_{s}$ is semisimple and $\lieg_r$ is the radical of $\lieg$.
  Then $\lieg_r$ and $\lieg$ are $\theta$-stable. Let $H$ be the
  Zariski closure of $G$ in $U\c$ and let $M$ be the Zariski closure
  in $U\c$ of the Lie subgroup of $G$ corresponding to $\lieg_r$. Then
  $H$ and $M$ are $\theta$-stable, so they are reductive. Since
  $\lieg_r$ is solvable, so is $M$, so that $M$ is a torus. Since
  $\lieg_r$ is $G$-stable, $M$ is normal in $H$.  Thus $M$ lies in the
  center of $H^0$. Going back to $\lieg$, we see that $\lieg_r$ is the
  center of $\lieg$ and that we can choose $\lieg_{s}=[\lieg,\lieg]$.
  Then $\theta$ respects 
  the Levi decomposition of $\lieg$ and
  $\liep=(\liep\cap\lieg_{s})\oplus(\liep\cap\lieg_r)$. 
  Since $\liea$ is maximal Abelian, it contains $\liep\cap\lieg_r$, so $ \lie a =\lie
  a'\oplus(\lieg_r\cap\liep)$ where $\lie a'$ is maximal commutative
  in $\liep\cap\lieg_{s}$. Then $\Ad(K)\lie a'=\liep\cap\lieg_{s}$
  \cite[Ch.\ 5 \S 6]{Helgason}, so that $\Ad(K)(\lie a)=\liep$.  It
  follows that $G=KAK$.
  \end{proof}

\begin{example} The Zariski closure of $G$  can be much larger than $G$. 
  For example, let $U\c=(\C^*)^2$, $U=S^1\times S^1$ and $G=\R^*$
  where $\lieg\simeq \R\subset\R^2\simeq i\lieu$ has irrational slope.
  Then the Zariski closure of $G$ in $U\c$ is $U\c$. The real
  algebraic closure of $G$ in $(\R^*)^2$ is $(\R^*)^2$.
\end{example}

Let $A$ and $\liea$ be as in Proposition \ref{KAK result}.  Then the
results in section \ref{commutative groups} show that every point in
$\mip\subset\m_{i\liea}$ has a neighborhood basis of
$\mu_{i\liea}$-adapted open sets.

The following lemma was established during discussions with H.
St\"otzel.

\begin{lemma} Let $\Omega_1$ and $\Omega_2$ be 
  $\mu_{i\lie a}$-adapted with $\Omega_1\cap\Omega_2=\emptyset$, and
  suppose that $W_i$ is a $K$-invariant subset of $\Omega_i$, $i=1$,
  $2$. Then $G\cdot W_1\cap G\cdot W_2=\emptyset$.
\end{lemma}

\begin{proof} It is enough to show that  
  $W_1\cap G\cdot W_2=\emptyset$. So assume that we have $g=k(\exp\xi)
  k'\in G=KAK$ and $w_i\in W_i$ such that $w_1=g\cdot w_2$. Then
  replacing $w_1$ by $k\inv w_1$ and $w_2$ by $k'w_2$, we can assume
  that $w_1=(\exp\xi)w_2$. Then $w_2\not\in\Omega_1$ and
  $(\exp\xi)w_2\in\Omega_1$, so there is a $t_1$ where $(\exp t\xi
  )w_2$ first enters $\Omega_1$. Similarly, $w_2\in \Omega_2$ while
  $(\exp \xi)w_2\not\in\Omega_2$, so that there is a last $t_2$ such
  that $(\exp t\xi)w_2\in\Omega_2$. By $\mu_{i\lie a}$-adaptedness,
  $t_1<t_2$, so that $(\exp t\xi)w_2\in\Omega_1\cap\Omega_2=\emptyset$
  for $t\in[t_1,t_2]$, a contradiction.
\end{proof}

\begin{corollary} \label{disjoint neighborhoods} Let 
  $C_1$ and $C_2$ be disjoint closed $K$-stable subsets of $\mip$.
  Then there are $G$-stable disjoint open subsets $\Omega_i\supset
  C_i$, $i=1$, $2$.  Hence the separation property holds.
\end{corollary}
\begin{proof} There are disjoint open subsets  
  containing the $C_i$, hence disjoint open $\mu_{i\lie a}$-adapted
  subsets $W_i$ containing the $C_i$, $i=1$, $2$. These in turn
  contain $K$-stable open neighborhoods $W_i'$ of $C_i$. Set
  $\Omega_1=G\cdot W_1'$ and $\Omega_2=G\cdot W_2'$. Then the
  $\Omega_i$ have the required properties.
\end{proof}

Now that we have the separation property, we have all the results of
\S \ref{semistable}. For completeness, we restate them here.

\begin{theorem} \label{first main theorem} Let $Z$ be a 
  holomorphic $U\c$-space with $U$-invariant K\"ahler form and moment
  mapping $\mu$, and let $G$ be a closed compatible subgroup of $U\c$.
  Then
\begin{enumerate}
\item The topological Hilbert quotient $\pi$ of $\mathcal{S}_G(\mip)$
  by $G$ exists, and the inclusion $\mip\to\mathcal{S}_G(\mip)$
  induces a homeomorphism $\mip/K\tosim \quot{\mathcal{S}_G(\mip)}G$.
  In particular, $ \quot{\mathcal{S}_G(\mip)}G$ is 
metrizable and locally compact.
\item If $A\subset \mathcal{S}_G(\mip)$ is closed, then
  $\pi(A)=\pi(A\cap \mip)$ is closed.
\item If $A_1$ and $A_2$ are closed and $G$-stable subsets of
  $\mathcal{S}_G(\mip)$, then $\pi(A_1)\cap\pi(A_2)=\pi(A_1\cap A_2)$.
\item If $G$ is commutative, then $\mathcal{S}_G(\mip)$ is open in $Z$
  and every $K$-orbit in $\mip$ has a neighborhood basis of
  $\muip$-adapted open sets.
\end{enumerate}
\end{theorem}
 
In case
$Z=\mathcal{S}_{U\c}(\m)$ we can establish more than the separation
property.

\begin{proposition}  Assume that the analytic Hilbert quotient 
  $\pi\colon Z\to \quot ZU\c$ exists and 
 that we have a $U$-invariant
  strictly plurisubharmonic function $\rho$ on $Z$, bounded below,
  such that $\pi\times \rho$ is proper. Let $\mu$ be associated with
  $\rho$ and let $x_0\in\mip$. Then $K\cdot x_0$ has a basis of
  $\muip$-adapted neighborhoods.  The elements of the neighborhood
  basis can be taken to be of the form $\Delta_{r_n}(\rho)\cap V_n$
  where $r_n>\rho(x_0):=r_0$, $\Delta_{r_n}(\rho)=\{z\in Z;\ 
  \rho(z)<r_n\}$ and $V_n$ is a $G$-invariant neighborhood of $x_0$.
\end{proposition}

\begin{proof} Fix a $K$-neighborhood $W_0$ of $K\cdot x_0$. Let 
  $W_0\supset W_1\supset W_2\supset \dotsb\supset K\cdot x_0$ be a
  basis of $K$-neighborhoods of $K\cdot x_0$ and let
  $r_1>r_2>\dotsb>r_0$ be such that $\lim r_n=r_0$. We show that
  $\Delta_{r_n}(\rho)\cap G\cdot W_n\subset W_0$ for $n$ sufficiently
  large.  By Proposition \ref{adapted sets}, the sets
  $\Delta_{r_n}(\rho)$ are $\mu$-adapted, hence the sets
  $\Delta_{r_n}(\rho)\cap G\cdot W_n$ are $\muip$-adapted and we have
  the proposition.
  
  Suppose that the proposition is false. Then we can assume that there
  are $x_n\in W_n$, $g_n\in G$ such that $g_n\cdot x_n\not\in W_0$
  while $\rho(g_n\cdot x_n)\leq r_n$. We may assume that $\lim
  x_n=x_0$. Since $\pi\times \rho$ is proper, we may assume that $\lim
  g_n\cdot x_n=y_0\in Z$. Let $z_0\in\overline{G\cdot y_0}\cap\mip$.
  We claim that $K\cdot z_0=K\cdot x_0$. If not, then by Corollary
  \ref{disjoint neighborhoods} we can find disjoint $G$-neighborhoods
  $\Omega_1$ and $\Omega_2$ of $z_0$ and $x_0$, respectively.
But we must have that
  $y_0\in \Omega_1$, so that $g_n \cdot x_n\in\Omega_1$ for $n$ large.
  Thus $x_n\in\Omega_1$ where $x_n\to x_0\in\Omega_2$. This is a
  contradiction.
  
  Now $\overline{G\cdot y_0}\cap\mip$ consists of the minima of $\rho$
  on the closed $G$-orbits in $\overline{G\cdot y_0}$, and we have just
  shown that these minima consist of the $K$-orbit of $x_0$. Since
  $\rho(x_0)=r_0$ and $\rho(y_0)\leq r_0$ we must have that $K\cdot
  y_0=K\cdot x_0\subset W_0$. But $g_n\cdot x_n\not\in W_0$ and
  $g_n\cdot x_n\to y_0$, a contradiction.
\end{proof}

\begin{theorem} \label{main theorem} Suppose that 
  $Z=\mathcal{S}_{U\c}(\m)$.  Then every point $z\in\mip$ has a
  neighborhood basis of $\muip$-adapted open subsets.
\end{theorem}

\begin{proof} By Theorem \ref{local plurisubharmonic} and the Proposition 
  above, every point in $\mip$ has a neighborhood basis of
  $\muip$-adapted open subsets.
\end{proof}

\section{Slices} 

In this section we assume that $Z$ is a holomorphic $U\c$-space with a
$U$-invariant K\"ahler structure and moment mapping $\mu$. We show
that we can find slices at points of $\mip$.  In particular, if the
analytic Hilbert quotient $\quot Z{U\c}$ exists (e.g., $Z$ is Stein),
then there are slices at points on closed $G$-orbits (see
Remark~\ref{slices closed orbits} below).  First we give a sufficient
condition that $K$-equivariant maps extend to be $G$-equivariant maps.

Let $G=K\exp\liep$ be a closed compatible subgroup of $U\c$.  A
$K$-stable subset $\Omega$ of $Z$ is said to be \emph{orbit convex
  with respect to $G$} if for every $\xi\in i\liep$ and $z\in\Omega$
the set $I_\Omega(z;\xi)=\{t\in\R;\ \exp it\xi\cdot z\in \Omega\}$ is
connected.  Note that a $\muip$-adapted subset of $Z$ is orbit convex
with respect to $G$.

Let $\Omega\subset Z$ be $K$-stable and open and let $X$ be a
topological $G$-space. Let $\varphi\colon\Omega\to X$ be
$K$-equivariant. We say that $\varphi$ is \emph{locally
  $G$-equivariant} if for every $x\in\Omega$ there is a neighborhood
$V_x$ of $e\in\exp(\liep)$ with $V_x\cdot x\subset\Omega$ such that
$\varphi(g\cdot x)=g\cdot\varphi(x)$ for $g\in V_x$.  

\begin{proposition} \label{orbit convex} 
  Let $A$ be a subgroup in $\exp\liep$ such that $G=KAK$.  Let
  $\Omega\subset\Omega_1\subset \Omega_2\subset Z$ be open sets where
  $\Omega$ and $\Omega_2$ are $K$-invariant and $\Omega_1$ is orbit
  convex with respect to $A$.  Let $X$ be a topological $G$-space and
  let $\varphi\colon \Omega_2\to X$ be a locally $G$-equivariant
  continuous map. Then there is a unique $G$-equivariant continuous
  map $\Phi\colon G\cdot \Omega\to X$ such that $\Phi|\Omega=\varphi$.
\end{proposition}

\begin{proof} For $z\in\Omega$ and $g\in G$, we define 
  $\Phi(g\cdot z)$ to be $g\cdot\varphi(z)$. This is clearly the
  desired mapping, as long as we can prove that it is well-defined. So
  let $z_1$, $z_2\in \Omega$ and $g_1$, $g_2\in G$ and suppose that
  $g_1\cdot z_1=g_2\cdot z_2$.  We need to show that $g_1\cdot \varphi
  (z_1)=g_2\cdot \varphi(z_2)$.  Equivalently, we have to show that
  $g\cdot\varphi(z_1)=\varphi(z_2)$ where $g=g_2\inv g_1$. Write
  $g=k(\exp\xi)k'$ where $k$, $k'\in K$ and $\xi\in\lie a$. Then
  $z_1':=k'\cdot z_1$ and $\exp\xi\cdot z_1'=z_2':=k\inv\cdot z_2$ are
  in $\Omega\subset\Omega_1$, hence $I_{\Omega_1}(z_1';-i\xi)$
  contains the interval $[0,1]$.  It follows that
  $\varphi(\exp(\xi)\cdot z_1')=\exp(\xi)\cdot\varphi(z_1')$ and from
  $K$-equivariance of $\varphi$ we finally obtain that $g\cdot
  \varphi(z_1)=\varphi(g\cdot z_1)=\varphi(z_2)$.
\end{proof} 

\begin{remark} \label{orbit convex remark} Suppose that $X$ is an analytic 
  $G$-space and that $\varphi$ is
  $\mathcal{C}^k$, $k\geq 1$, or real analytic. Then $\Phi$ is
  $\mathcal{C}^k$ or real analytic. If $X$ is a complex $G$-space and
  $\varphi$ is holomorphic, then $\Phi$ is holomorphic.
\end{remark}

The following is an analogue in our setting of Luna's fundamental
lemma.  

\begin{lemma}\label{fundamental lemma} 
  Let $X$ be a topological $G$-space and let $x\in X$.
  Let $\varphi\colon X\to Z$ be a $G$-equivariant continuous mapping
  which maps $K\cdot x$ injectively into $G\cdot z$ where
  $z=\varphi(x)$.  Assume that $\varphi$ is a local homeomorphism at
  $x$ and that $z\in\mip$. Then there is an open $G$-stable
  neighborhood $W$ of $G\cdot x$ which is mapped homeomorphically by
  $\varphi$ onto the open $G$-neighborhood $\varphi(W)$ of $G\cdot z$.
\end{lemma}
 
\begin{proof} By equivariance,  $\varphi$ 
  is a local homeomorphism along $K\cdot x$. It follows that there is
  a neighborhood $\Omega$ of $K\cdot x$ which is mapped
  homeomorphically onto its image $\varphi(\Omega)\subset Z$.  Let
  $\Omega_1$ be a $K$-stable neighborhood of $x$ which is contained in
  $\Omega$. Then $\varphi | \Omega_1$ is a locally $G$-equivariant
  homeomorphism.  Let $A\subset\exp\liep$ be a subgroup such that
  $G=KAK$. Then by Proposition~\ref{mip-adapted sets commutative case}
  and Lemma~\ref{mip union lemma} we may find an orbit convex (with
  respect to $A$) neighborhood $V_1$ of $K\cdot z$ inside $V_2:=
  \varphi( \Omega_1)$. Let $V$ be a $K$-invariant neighborhood of
  $K\cdot z$ which is contained in $V_1$.  Since $(\varphi |
  \Omega_1)\inv | V_2$ is a locally $G$-equivariant mapping, there is
  a unique $G$-equivariant extension $\psi$ of $(\varphi |
  \Omega_1)\inv | V$ to $G\cdot V$. The composition $\varphi\circ
  \psi$ is the identity on $G\cdot V$ by construction. Similarly,
  $\psi\circ\varphi$ is the identity on $W:=G\cdot\psi(V)\subset
  G\cdot\Omega_1$.
\end{proof}

\begin{remark} \label{differentiability remark} Suppose that $X$ is a real analytic (resp.\ complex) $G$-space. If
$\varphi$ is a local analytic diffeomorphism (resp.\ 
  local biholomorphism), then $\varphi | W$ is an analytic
  diffeomorphism (resp.\ biholomorphism). If $\varphi$ is a local
  $\mathcal{C}^k$ diffeomorphism, $k\geq 1$, then $\varphi | W$ is a
  $\mathcal{C}^k$ diffeomorphism.
\end{remark}

\begin{remark} \label{fundamental lemma remarks} If $G=U\c$ and
  $\varphi$ is holomorphic, then Lemma~\ref{fundamental lemma} is a
  strict generalization of the usual fundamental lemma in the complex
  analytic category. Similarly, Theorem~\ref{smooth case slice},
  Theorem~\ref{slice theorem} and Proposition~\ref{equivariant
    isomorphism} below give statements in the complex analytic
  category. Proposition~\ref{equivariant isomorphism} in this setting
  may even be new.
\end{remark}

We give a proof of the following well-known result.

\begin{lemma}\label{fixed point slice theorem} Let  $z\in Z$ be 
  a $U\c$-fixed point.  Then there is a $U\c$-stable neighborhood $W$
  of $z$ and a $U\c$-equivariant biholomorphic map $\psi\colon
  W\to\tilde Z$ where $\tilde Z\subset T_z(Z)$ is a $U\c$-stable
  closed subspace of a $U\c $-stable open subset of $T_zZ$.
\end{lemma}

\begin{proof} Let $\psi$ be a $U$-equivariant biholomorphic map of a 
  $U$-neighborhood $W'$ of $z\in Z$ onto a locally closed $U$-stable
  subspace $Z'$ of $T_zZ$ where $\psi(z)=0$.  We may assume that $Z'$
  is closed in the ball $B$ of radius 1 for some $U$-invariant norm on
  $T_zZ$. Then $B$ is $\mu'$-adapted for $\mu'$ the moment mapping
  associated to the square of the norm function on $T_zZ$.  Moreover,
  $Z'$ is locally $U\c$-stable (see Remark~\ref{locally g stable}). It
  follows that $\tilde Z:=U\c\cdot Z'$ is a closed analytic
  $U\c$-subspace of the $U\c$-stable open subset $U\c\cdot B$ such
  that $\tilde Z\cap B=Z'$.  Note that $\psi\colon W'\to \tilde Z$ is
  automatically locally $U\c$-equivariant. Changing $\mu$ by a
  constant, we may assume that $z\in\m$. Then $z$ has a neighborhood
  basis of $\mu$-adapted open sets (Theorem~\ref{main theorem}) and by
  Proposition~\ref{orbit convex} we may assume that $\psi$ is defined
  on a $U\c$-invariant neighborhood of $z$. Now apply
  Lemma~\ref{fundamental lemma}.
\end{proof}

We now need a result on complete reducibility of certain
representations of $G$.

\begin{lemma} \label{completely reducible lemma} Let $V$ be a 
  complex representation space of $U\c$. Then $V$ is completely
  reducible as a real representation of $G$.
\end{lemma}

\begin{proof}  We may assume that $V\simeq\C^n$ has the usual hermitian inner product   $\langle\, ,\rangle$ and that $U$ is
  a subgroup of $U(n,\C)$. The real part of $\langle\, ,\rangle$ is
  the usual real inner product $(\, ,)$ on the underlying real vector
  space $\R^{2n}$. Then $K$ consists of orthogonal matrices and
  elements of $\exp(\liep)$ are real symmetric matrices.  Let $W$ be a
  real $G$-submodule of $V$.  Taking perpendicular relative to $(\,
  ,)$ we have the $K$-stable subspace $W^\perp$.  Moreover, if $w\in
  W$, $w^\perp\in W^\perp$ and $g\in\exp(\liep)$, then $0=(g\cdot
  w,w^\perp)=(w,g\cdot w^\perp)$, so that $W^\perp$ is also
  $\exp(\liep)$-stable. Thus $W^\perp$ is $G$-stable.
\end{proof}

\begin{remark} \label{singular kaehler remark} Let $z\in Z$ and let $H$ be a compatible reductive
subgroup of
$U\c$ fixing
$z$. Then locally and $H$-equivariantly we may consider $Z$ as an $H$-stable subset of $T_zZ$.
The $U$-invariant K\"ahler structure on $Z$ extends locally to an $(H\cap U)$-invariant K\"ahler
structure on $T_zZ$ \cite{Narasimhan}.

\end{remark}
 
\begin{corollary}\label{Nx corollary} Let $z\in Z$ such that 
  $G_z$ is a compatible subgroup of $U\c$ (e.g., $z\in\mip$). Then the
  representation of $G_z$ on $T_z Z$ is completely reducible. In
  particular, there is a direct sum decomposition $T_z Z=T_z(G\cdot
  z)\oplus N_z$ where $N_z$ is $G_z$-stable.
\end{corollary}  
\begin{proof} 
  The Zariski closure of $G_z$ in $U\c$ is  reductive (Remark \ref{reductive closure remark}),
and lies in $(U\c)_z$. Now applying the argument of Lemma \ref{completely reducible lemma} we
can choose $N_z=(\lieg\cdot x)^\perp$ for the inner product on $T_zZ$ coming from the K\"ahler
structure.
\end{proof}

\begin{theorem} \label{smooth case slice} Assume that $Z$ 
  is smooth at $z\in\mip$. Then there is a geometric $G$-slice at $z$,
  i.e., there is a $G_z$-invariant locally closed real analytic
  submanifold $S$ of $Z$, $z\in S$, such that $G\cdot S$ is open in
  $Z$ and such that the natural map $\twist G{G_z}S\to G\cdot S$,
  $[g,s]\mapsto g\cdot s$, is a real analytic isomorphism.
\end{theorem}

\begin{proof} 
  Let $H$ denote the analytic (hence algebraic) Zariski closure of
  $G_z$ in $U\c$.  Note that $z$ is an $H$-fixed point. By the
  holomorphic slice theorem for fixed points (Lemma~\ref{fixed point
    slice theorem}) there is an $H$-stable open neighborhood $D\subset
  T_z(Z)$ and an $H$-equivariant open embedding $\imath \colon D\to Z$
  with $\imath(0)=z$. From Corollary~\ref{Nx corollary} we know that
  there is a $G_z$-complement $N_z$ to $T_z(G\cdot z)$ in $T_z(Z)$.
  Let $\hat S:=\imath(D\cap N_z)$ and let $\varphi\colon \twist
  G{G_z}{\hat S}\to Z$ denote the natural map.  Then $\varphi$ is
  $G$-equivariant, a local diffeomorphism at $[e,z]$ and an
  isomorphism on the $G$-orbit through $[e,z]$.  Applying
  Lemma~\ref{fundamental lemma} we obtain a $G$-neighborhood $\Omega$
  of $[e,z]\in\twist G{G_z}{\hat S}$ which is mapped isomorphically
  onto $G\cdot \varphi(\Omega)$. But $\Omega=\twist G{G_z}S$ where
  $S=\Omega\cap \hat S$.
\end{proof}

\begin{remark} 
By Lemma \ref{elementary}(2) we have $(\lie g\cdot z)^\perp=(\lie k\cdot z)^\perp\cap
\ker(d\muip(z))$. Thus, on the
  tangent space level, the slice for the $G$-action is the
  intersection of $\mip$ with a slice for the $K$-action.
\end{remark}

The slice theorem for the case that $z$ is not a smooth point of $Z$
requires more work.  If we knew that $Z$ was a $G$-subspace of a
smooth $G$-space $\hat Z$, then the slice theorem for $\hat Z$ would
imply the slice theorem for $Z$. But the only way we see to find a
$\hat Z$ is to prove the slice theorem for $Z$!

Assume that $z\in\mip$ and let $H$ denote the Zariski closure of $G_z$
in $U\c$.  By Lemma~\ref{fixed point slice theorem} we have an
$H$-equivariant embedding $\psi$ of an open neighborhood $W$ of $z$
onto a locally analytic subset $\tilde Z\subset T_zZ$ where
$\psi(z)=0$. Write $T_zZ=T_z(U\c\cdot z)\oplus N$ where $N$ is
$H$-stable.
The set $\tilde \Sigma:=\tilde Z\cap N$ is an $H$-stable locally
analytic subset of $N$ and $\Sigma:=\psi\inv(\tilde \Sigma)$ is an
$H$-stable analytic subset of $W$.

We have an $\Ad H$-stable
decomposition $\lieu\c=(\lieu\c)_z\oplus\liem$.  Note that $\liem \to
U\c/(U\c)_z$, $m\mapsto \exp(m)(U\c)_z$, is biholomorphic on an open
neighborhood of $0\in\liem$.

\begin{lemma} \label{S0 lemma}
  The holomorphic map $\varphi^0\colon \liem\times \Sigma\to Z$,
  $(m,s)\mapsto \exp(m)\cdot s$, is $H$-equivariant and there is an
  open $H$-stable neighborhood $V$ of $(0,z)\in\liem\times \Sigma$
  such that $\varphi^0$ maps $V$ biholomorphically onto an open subset
  of $Z$.
\end{lemma}

\begin{proof} 
  By construction, $d\varphi^0_{(0,z)}$ is injective. This implies
  that $\varphi^0$ maps an open neighborhood $V=\Omega^0\times
  \Sigma^0$ of $(0,z)$ biholomorphically onto a closed analytic subset
  of a neighborhood $Z^0$ of $z\in Z$. Here $\Omega^0$ is a connected
  neighborhood of $0\in\liem$ and $\Sigma^0$ is a neighborhood of
  $z\in\Sigma$. Let $Z_\beta^0$, $\beta\in B$, be the irreducible
  components of $Z^0$.  We may assume that each $Z_\beta^0$ is
  irreducible at $z$. There is an open neighborhood $\Omega$ of $e\in
  U\c$ such that $\Omega\cdot z$ is a locally closed submanifold of
  $Z_\beta^0$ for each $\beta$.  Set $\tilde Z^0:=\psi(Z^0)$ and
  $\tilde Z^0_\beta:=\psi(Z^0_\beta)$, $\beta\in B$.  For each
  $\beta$, $\psi(\Omega\cdot z)\subset \tilde Z_\beta^0$ is a smooth
  submanifold through $0$ which is transversal to $N$, so we have
\[
\dim \tilde \Sigma_{\alpha\beta}^0= \dim Z_\beta^0- \dim \Omega\cdot
z=\dim Z_\beta^0-\dim \liem
\] 
for every irreducible component $\tilde\Sigma_{\alpha\beta}^0$ of
$\tilde Z_\beta^0\cap N$ which contains $0$.  For each $Z_\beta^0$
choose an irreducible component $\Omega^0\times
\Sigma_{\alpha\beta}^0$ of $V$ such that $z\in
\Sigma_{\alpha\beta}^0\subset Z_\beta^0$. Since $\dim (\liem\times
\Sigma^0_{\beta\alpha})=\dim Z_\beta^0$ we have
$\varphi(\Omega^0\times \Sigma_{\alpha\beta}^0)=Z^0_\beta$.  Hence
$\varphi^0$ is a surjective closed embedding of $\Omega^0\times
\Sigma^0$ onto $Z^0$ and therefore biholomorphic.  
\end{proof}
 \begin{remark} We have $T_z\Sigma=N\subset T_zZ$.
\end{remark}

Since $G_z$ is a compatible subgroup of $U\c$, the action of $G_z$ on
$\lieu\c$ is completely reducible. Let $\liem'\subset\liem$ be an
$(\Ad G_z)$-stable complement to $(\lieu\c)_z+\lieg$ in $\lieu\c$.  We
have the twisted product $\twist G{G_z}(\liem'\times \Sigma)$ where
$G_z$ has the product action on $\liem'\times \Sigma$. There is the
$G$-equivariant map $\varphi\colon \twist G{G_z}(\liem'\times
\Sigma)\to Z$, $(g,m,\sigma)\mapsto g\exp(m)\cdot\sigma$.

\begin{proposition} \label{slice proposition} 
  There is a  $G_z$-invariant open neighborhood $S$ of
  $(0,z)\in\liem'\times\Sigma$ such that $\varphi\colon\twist G{G_z}
  S\to Z$ is a $G$-equivariant open embedding.
\end{proposition}

\begin{proof} By  Lemma~\ref{fundamental lemma} and Remark 
  \ref{differentiability remark} it is enough to 
  prove that $\varphi$ is an isomorphism in a neighborhood of
  $[e,(0,z)]$. Let $\liem''$ be an $(\Ad G_z)$-stable complement to
  $\liem'$ in $\liem$. Note that $\liem''\simeq\lieg/\lieg_z$.  It is
  enough to show that
  $$
  \varphi'\colon\liem''\oplus\liem'\times\Sigma\to Z, \ 
  (m'',m',\sigma)\mapsto \exp(m'')\exp(m')\cdot\sigma
  $$
  is a diffeomorphism in a neighborhood of $(0,z)$. Since
  $\varphi^0$ is biholomorphic, it is enough to show that
  $\sigma:=(\varphi^0)\inv\circ\varphi'\colon\liem\times\Sigma\to\liem\times\Sigma$
  is a diffeomorphism near $(0,z)$.  Clearly $d\sigma_{(0,z)}$ is the
  identity. Now we just apply the lemma below.
\end{proof}

\begin{lemma} 
  Let $X$ be a germ at $0\in \C^n$ of a complex analytic set. Let
  $\sigma\colon X\to X$, $\sigma(0)=0$, be a germ of a real analytic
  mapping such that $\sigma(X)\subset X$. Assume that $d\sigma_0\colon
  T_0X\to T_0X$ is an isomorphism. Then $\sigma(X)=X$ (as germs).
\end{lemma}

\begin{proof} We may assume that $T_0X=\C^n$. Then we may assume that 
  a representative of $\sigma$ extends to a real analytic
  diffeomorphism $\tau\colon W_1\to W_2$ where $W_1$ and $W_2$ are
  neighborhoods of $0$. We may assume that $X$ is represented by
  closed analytic subsets $X_j$ of $W_j$, $j=1$, $2$, such that $\tau(
  X_1)\subset X_2$.
  
  We first assume that the $X_j$ are irreducible. Then each
  $X_j\setminus \Sing X_j$ is connected and smooth where $\Sing X_j$
  has real codimension at least two in $X_j$. There is certainly a
  smooth point $x_1\in X_1$ such that $x_2:=\sigma(x_1)$ is a smooth
  point in $X_2$. Since $\tau$ is an isomorphism, $\sigma$ is open at
  $x_1$. Hence there is a neighborhood $V_2$ of $x_2$ in
  $X_2\setminus\Sing X_2$ such that $\tau\inv(V_2)\subset X_1$. By the
  identity principle for real analytic maps we have
  $\tau\inv(X_2)\subset X_1$. Hence $\sigma$ is onto.
  
  In general, each $X_j\setminus \Sing X_j$ is a finite union of $d_j$
  connected components whose closures are the irreducible components
  of $X_j$, where $d_1=d_2$.  By the argument above, if $Y$ is an
  irreducible component of $X_1$, then $\sigma(Y)$ is an irreducible
  component of $X_2$. Since $\tau$ is invertible, the irreducible
  components of $X_1$ are mapped onto the irreducible components of
  $X_2$.  Hence $\sigma(X_1)= X_2$.
\end{proof}

Proposition \ref{slice proposition} has several important
corollaries.
\begin{corollary}\label{closed orbits of G} Let $z\in\mip$. Then $G\cdot z$ is  closed in
  $\mathcal{S}_G(\mip)$ and closed in a $G$-stable open subset of $Z$.
  In particular, a $G$-orbit in $\mathcal{S}_G(\mip)$ is closed if and
  only if it intersects $\mip$.
\end{corollary}
\begin{proof} By the proposition, $G\cdot z$ is closed in $G\cdot S$. 
  Let $y\in\overline{G\cdot z}\setminus G\cdot z$.  If
  $y\in\mathcal{S}_G(\mip)$, then $\overline{G\cdot y}\cap\mip$ must
  be $K\cdot z$ by Corollary \ref{intersection with mip corollary},
  but this is impossible since $y$ lies in the complement of $G\cdot S $.
\end{proof}

\begin{corollary} For every $z\in \mip$, the saturation $\mathcal{S}_G(\{z\})$
  is closed in an  
open $G$-stable neighborhood of
  $z$.
\end{corollary} 

\begin{proof} We may $G_z$-equivariantly identify 
$S\subset\liem'\times\Sigma$ with its image in $\liem'\times N\subset T_zZ$. Then $G_z$ is
compatible with the unitary structure on $T_zZ$ and we have a topological Hilbert quotient
$\pi\colon T_zZ\to\quot{T_zZ}{G_z}$.  The null cone $\tilde N:=\pi\inv(\pi(0))$ consists of the 
$y\in T_zZ$
  such that the closure of $G_z\cdot y$ contains $0$. Then
$\tilde N\cap S$ is closed in $S$, and $\mathcal{S}_G(\{z\})\cap G\cdot
  S\simeq \twist G{G_z}(\tilde N\cap S)$ is closed in $G\cdot S$.
\end{proof}

Similarly, we have

\begin{corollary}\label{minimal dimension orbit}Let $y\in\mathcal{S}_G(\mip)$
  and let $z\in\overline{G\cdot y}\cap\mip$. Then $G\cdot z$ is the
  unique orbit of minimal dimension in $\overline{G\cdot y}\cap
  \mathcal{S}_G(\mip)$.
\end{corollary}

\begin{corollary} 
  Let $z\in\mip$.  Then there is a smooth real analytic $G$-space
  $\hat Z$ and a $G$-equivariant embedding of a $G$-neighborhood of
  $z$ into $\hat Z$.
\end{corollary}

\begin{proof} 
  We have an embedding of a $G$-neighborhood of $z$ into $\twist
  G{G_z}(\liem'\times N)$.
\end{proof}

\begin{corollary} Let $z\in\mip$. Then there is a $G$-invariant 
  neighborhood $W$ of $z$ such that the topological Hilbert
  quotient $\quot WG$ exists.
\end{corollary}
 
\begin{proof} As before, 
the topological Hilbert quotient
$\pi\colon  T_zZ\to\quot{T_zZ}{G_z}$   exists. 
We may identify $S$ with its image in $T_zZ$. 
Let $S'$ denote
$S\cap\pi\inv(\pi(\overline{S})\setminus\pi(\overline{S}\setminus S))$. Then $0\in S'$, $S'$ is $G_z$-stable and open in $S$, 
and
$S'\subset\mathcal{S}_G(\mip\cap S')$. It follows that the topological Hilbert quotient of $S'$ exists, hence so does the
topological  Hilbert quotient of
$W:=\twist G{G_z}S'$.

\end{proof}

\begin{corollary} \label{semianalytic corollary} 
 Assume that $Z':=\mathcal{S}_G(\mip)$ is open and
 let $z\in\mip$. 
  Then a neighborhood of $G\cdot z\in\quot {Z'}G$ is
  homeomorphic to a semianalytic set.
\end{corollary}

\begin{proof}  
  We continue from the proof above. We may assume that $S\subset Z'$. 
Then $\quot{(\twist
G{G_z}S)}G$ parametrizes the closed $G$-orbits, so we need only 
show that the closed $G_z$-orbits in $S$ are
locally parametrized by a semianalytic set. 
We may consider $S$ as a subset of $V:=T_zZ$.
We have a unitary structure on $V$ for which the image of $G_z$ in
$\GL(V)$ is compatible. Let
  $p_1,\dotsc, p_d$ be generators of $\R[V]^{K_z}$ and let
  $p=(p_1,\dotsc,p_d)\colon V\to \R^d$. Then $V/K_z$ is homeomorphic to
  the closed semialgebraic subset $p(V)$. Moreover,  
$\mip\subset V$ is semialgebraic, hence $p(\mip)$ is also
  semialgebraic and parametrizes the closed $G_z$-orbits in $V$.  
We can find a nonnegative
real valued $K_z$-invariant real analytic
  function $f$ (near $0$) whose zeroes define $S$ near $0$. Then $f=p^*h$ where
  $h$ is analytic in a neighborhood of $0\in\R^d$. The zeroes of $h$
  on $p(\mip)$ give the germ of the image of $\mip\cap S$, where $(\mip\cap S)/K_z$
parametrizes the closed $G_z$-orbits in $S$. Thus a
  neighborhood of $0\in\quot S{G_z}$ is homeomorphic to a semianalytic
  set.
\end{proof}

\begin{remark} \label{sgmip open remark} Let $Z$ be a
projective variety with a $U\c$-linearized very ample line bundle $L$. Let $V$ denote the dual of the  $U\c$-module $\Gamma(Z,L)$.
Then $Z$ embeds into $P(V)$ and $\mathcal{S}_G(\mip)$ is the  intersection of $Z$ with the image of $V\setminus
\mathcal{N}$ in $P(V)$, where $\mathcal{N}$ is the null cone of the $G$-action on $V$. Hence $\mathcal{S}_G(\mip)$ is open in
$Z$.
\end{remark}

We may reformulate Proposition~\ref{slice proposition} as follows.

\begin{theorem} \label{slice theorem} 
  For every $z\in\mip$ there is a locally closed real analytic
  $G_z$-stable subset $S$ of $Z$, $z\in S$, such that the natural map
  $\twist G{G_z}S\to Z$ is a real analytic $G$-isomorphism onto the
  open set $G\cdot S$.
\end{theorem}

\begin{remark} \label{saturated slice remark} 
  Assume that $Z':=\mathcal{S}_G(\mip)$ is open. Then we may assume that $G\cdot S\subset Z'$.
  Let $\pi \colon Z'\to \quot {Z'}G$ be the quotient mapping. Then the
  image $C$ of the complement of $G\cdot S$ in $Z'$ is closed and does not
  contain the closed orbit $G\cdot z$.  Thus replacing $G\cdot S$ by
  the inverse image of the complement of $C$ we may arrange that
  $G\cdot S$ is $G$-saturated.  It follows automatically that $S$ is
  $G_z$-saturated.
\end{remark}
 
\begin{example} This is a continuation of Example~\ref{Lie algebra example}.
  Let $\pi\colon \lieg \to \quot\lieg G$ denote the quotient map. The
  slice theorem implies that $\pi\inv(\pi(x))=G\cdot x$ for every
  $x\in X$.  In particular, the orbit space $X/G\cong \pi(X)$ is
  Hausdorff and the quotient map is given by restricting $\pi$ to $X$.
\end{example}

We now can state some variants of Theorem~\ref{compact slice theorem} and 
Proposition~\ref{real proper}.
For a general holomorphic $U\c$-space $Z$, let $\Comp_G(Z)$ denote
$\{z\in Z; \text{ for all $w\in \overline{G\cdot z}$ the isotropy
  group $G_w$ is compact}\}$. 

\begin{proposition} \label{compactorbits2}
  If $Z':=\mathcal{S}_G(\mip)$ is open in $Z$, then $\Comp_G(Z')$ is
  open in $Z'$ and the natural map 
$\twist GK(\mip\cap\Comp_G(Z'))\to\Comp_G(Z')$ is a homeomorphism. In particular,
the $G$-action on $\Comp_G(Z')$ is proper. 
\end{proposition}
\begin{proof} Let $C$ denote $\Comp_G(Z')$. It follows
from the slice theorem that $C$ is 
open in $Z'$ and that every $G$-orbit in $C$ is closed.  Hence $\twist
GK(\mip\cap C )\to C$ is one to one and onto. We have a homeomorphism
$(\mip\cap C)/K\to\quot{C}G$. Let $\pi\colon C\to \quot CG$ denote the 
quotient mapping.  Suppose that
$[g_n,z_n]\in\twist GK(\mip\cap C)$ such that $g_n\cdot z_n$ converges in 
$C$. Then $\pi(z_n)$ is convergent, so  we can
assume that $z_n\to z\in\mip\cap C$. If $S$ is a slice at $z$, then 
the action $G\times G\cdot S\to G\cdot S$ is proper. Hence
we can assume that $g_n\to g\in G$. Thus $\twist GK(\mip\cap C)\to C$ 
is a homeomorphism.
\end{proof}
\begin{remark} \label{slices closed orbits} Suppose that the analytic Hilbert 
  quotient $\quot Z{U\c}$ exists and let $G\cdot z$ be a closed orbit.
  Then we may assume that $Z=\mathcal{S}_{U\c}(\m)$ (see
  Remark~\ref{Hilbert quotient remark}), so that $G\cdot z$ intersects
  $\mip$.  Hence there is a   slice at $z$.
\end{remark}

>From Remark \ref{Hilbert quotient remark} and Theorem~\ref{slice
  theorem}   we   obtain

\begin{proposition} \label{compactorbits}
  If the analytic Hilbert quotient $Z\to \quot Z U\c$ exists, then
  $\Comp_G(Z)$ is open in $Z$ and the natural map $\twist
  GK(\mip\cap\Comp_G(Z))\to\Comp_G(Z)$ is a homeomorphism.  In
  particular, the $G$-action on $\Comp_G(Z)$ is proper.
\end{proposition}

Here is a criterion for an equivariant map to be a homeomorphism.

 \begin{proposition} \label{equivariant isomorphism} Let $X$ be a Hausdorff topological $G$-space
  such that the topological Hilbert quotient $\pi\colon X\to \quot X
  G$ exists and such that $\quot XG$ parametrizes the closed
  $G$-orbits.  Assume that $Z=\mathcal{S}_{G}(\mip(Z))$ and that we
  have a continuous mapping $\varphi\colon X\to Z$ and the following
  properties.
\begin{enumerate}
\item $\varphi$ is a local homeomorphism.
\item $\varphi$ is $G$-equivariant.
\item $\quot{\varphi}G$ is a bijection.
\item There is a $K$-orbit $K\cdot x_0$ which is mapped isomorphically
  onto $K\cdot z_0\subset \mip(Z)$.
\item $\quot XG$ is connected.
\end{enumerate}
Then $\varphi$ is a homeomorphism.
\end{proposition}

\begin{proof} Let $G\cdot z$ be a closed orbit
  in $Z$. By Corollary \ref{closed orbits of G} we may assume that
  $z\in\mip(Z)$. Since $\varphi$ is a local homeomorphism, the fiber
  $\varphi\inv(z)$ is discrete and consequently each $G$-orbit in
  $\var\inv(G\cdot z)\cong \twist G {G_z}\varphi\inv(z)$ is open in
  $\varphi\inv(G\cdot z)$. Thus each $G$-orbit in $\var\inv(G\cdot z)$
  is closed in $\var\inv(G\cdot z)$ and also closed in $X$.  But
  $\quot \varphi G$ is a bijection, so it follows that $G\cdot
  x=\varphi\inv(G\cdot z)$ for some $x\in \varphi\inv(z)$ and then we
  claim that $\var\inv(K\cdot z)=K\cdot x$. For suppose that
  $\varphi(k\cdot\exp(\xi)\cdot x)=k'\cdot z$ where $k$, $k'\in K$ and
  $\xi\in\liep$. Then $\exp(\xi)\cdot z\in K\cdot z$ which implies
  that $\xi_Z$ vanishes at $z$. Since $(\varphi\vert G\cdot x) \colon
  G\cdot x \to G\cdot z$ is automatically smooth and is a local
  diffeomorphism, it follows that $\xi_X$ vanishes at $x$, and we have
  the claim. Note that since $\quot\varphi G$ is a bijection, the
  closed $G$-orbits in $X$ are precisely the inverse images of the
  closed $G$-orbits in $Z$.
  
  We now show that $\quot \varphi G$ is a homeomorphism. Let $x\in X$
  such that $G\cdot x$ is closed, let $V$ be a neighborhood of $G\cdot
  x$ in $\quot XG$ and set $W:=\pi_X\inv(V)$.  Since $\varphi$ is open
  and $G\cdot\varphi(x)$ is closed, $\pi_Z(\varphi(W))$ contains a
  neighborhood of $G\cdot \varphi(x)$ (Remark \ref{open map remark}).
  Hence $\quot\varphi G$ is open at $G\cdot x$ and $\quot\varphi G$ is
  a homeomorphism.
  
  Now let $\mip(X)'$ denote the inverse image of $\mip(Z)$ in $X$.
  Then $\mip(X)'\to \mip(Z)$ is an open mapping, hence so is the
  composition $\mip(X)'\to\quot ZG$ and the quotient mapping
  $\mip(X)'\to\quot XG$.  It follows that the inclusion $\mip(X)'\to
  X$ induces a homeomorphism $\mip(X)'/K\tosim \quot XG$. Then
  $\varphi | \mip(X)'$ is proper and a local homeomorphism, hence a
  covering map. Now (3), (4) and (5) imply that the covering has one
  sheet, i.e., $\varphi | \mip(X)'$ is a homeomorphism onto $\mip(Z)$.
  Let $x\in\mip(X)'$.  Since $\varphi | (K\cdot x)\colon K\cdot x\to
  K\cdot \varphi(x)$ is an isomorphism, Lemma~\ref{fundamental lemma}
  shows that $\varphi$ is a homeomorphism on a $G$-neighborhood of
  $G\cdot x$.  Since $\quot\varphi G$ is a homeomorphism, it follows
  that $\varphi$ is a homeomorphism.
\end{proof}

\begin{remark} Via the homeomorphism $\varphi$ we can impose a complex
  $G$-space structure on $X$ such that $\varphi$ is biholomorphic.
\end{remark}
 
\section{The Hilbert Mumford Criterion} 
Let $G$ be a compatible subgroup of $U\c$ and choose a closed subgroup
$A\subset\exp\liep$ such that $G=KAK$.  Let $X$ be a $G$-space such
that the topological Hilbert quotient $\pi\colon X\to \quot X A$
exists and is regular. That is, $\quot XA$ is Hausdorff and if $\xi$
is a point of $\quot XA$ and $C$ is a closed subset not containing
$\xi$, then there are disjoint neighborhoods of $\xi$ and $C$. We also
assume that $\pi$ parametrizes the closed $A$-orbits and maps closed
$A$-stable subsets to closed subsets of $\quot X A$.

We use an argument of Richardson to establish the main part of the
Hilbert Mumford criterion.

\begin{proposition} \label{propositionrichardson}
  Let $X$ be as above and let $Y$ be a closed $G$-stable subset of
  $X$. Then
  \[
  \mathcal{S}_G(Y)=K\cdot \mathcal{S}_A(Y).
  \]
\end{proposition}

\begin{proof}
  Since $\mathcal{S}_G(Y)$ is $K$-stable and
  $\mathcal{S}_A(Y)\subset\mathcal{S}_G(Y)$, we have $K\cdot
  \mathcal{S}_A(Y)\subset \mathcal{S}_G(Y)$.  Now let $z\in
  \mathcal{S}_G(Y)$. If $\pi(K\cdot z)\cap\pi(Y)=\emptyset$, then it
  follows from regularity of $\quot XA$ that there are open
  $\pi$-saturated disjoint neighborhoods $\Omega_2$ of $Y$ and
  $\Omega_1$ of $K\cdot z$. Let $\Omega_1'$ be a $K$-stable
  neighborhood of $K\cdot z$ in $\Omega_1$ and define $\Omega_2'$
  similarly. Then $G\cdot\Omega_1'\cap G\cdot\Omega_2'=\emptyset$, so
  that $z\not\in\mathcal{S}_G(Y)$, which is a contradiction. Thus
  there is a $k\in K$ such that $\overline{A\cdot k\cdot z}\cap Y\ne
  \emptyset$, i.e., $z\in K\cdot \mathcal{S}_A(Y)$.
\end{proof}

Let $A$ be a commutative connected simply connected real Lie group
with Lie algebra $\lie a$. Then the exponential map gives an
isomorphism of $\lie a$ with $A$.  We consider here only finite
dimensional continuous (hence real analytic) representations of $A$. We say that a real  $A$-module $W$   is
{\it completely $1$-reducible\/} if it is completely reducible with each irreducible component being of dimension 1.  Then every
isotropy group of $W$  is connected (and therefore simply connected) and for
every $x\in W$, the $A_x$-module $W$ is completely 1-reducible.

\begin{remark} \label{Arepremark}
  Let $W$ be a real $A$-module.  If $W$ is completely 1-reducible, then
  one can choose an isomorphism $W\simeq \R^n$ such that the image
  $A'$ of $A$ in $\GL(n,\R)$ consists of real positive diagonal
  matrices.  Then $A'\subset\exp i\lieu(n,\C)$ is a compatible subgroup of $\GL(n,\C)$ with
  its usual Cartan decomposition. Conversely, given an isomorphism
  $W\simeq\R^n$ such that the image of $A$ in $\GL(n,\C)$ is compatible and lies in $\exp i\lieu(n,\C)$, then the image
consists of symmetric real matrices, so
  that $W$ is completely 1-reducible.
\end{remark}

\begin{example}
  Let $W$ be completely 1-reducible with isotypic decomposition
  $W=W_{\chi_0}\oplus\dots\oplus W_{\chi_r}$ where $\chi_0$ is the
  trivial character. Then $W_{\chi_0}=W^A$.  We may assume that $W$ is
  a real subspace of a holomorphic $U\c$-representation such that the
  $U$-representation is unitary with respect to some Hermitian inner
  product $\langle\ , \ \rangle$ with associated norm $\lVert\ 
  \rVert$. We may assume that $A\subset \exp(i\lieu)$.  Let
  $\mu_{i\lie a}$ denote the restriction to $W$ of the $i\lie
  a$-component of the moment map associated with the strictly
  plurisubharmonic $U$-invariant exhaustion function
  $\rho:=\frac{1}{2}\lVert\ \rVert^2$.  A simple calculation shows
  that
  \[
  \mu^\xi(z)=i\langle \xi z, z\rangle\quad \xi\in \lieu, z\in W.
  \]
  For $z=z_0+\dotsb +z_r\in W_{\chi_0} \oplus\dotsb\oplus W_{\chi_r}$
  we obtain
  \[
  \mu_{i\lie a}(z)=i(\lVert z_1\rVert^2\chi_1+ \dotsb +\lVert
  z_r\rVert^2\chi_r).
  \]
 Hence
  \[C:=\mu_{i\lie a}(W)=\{i(a_1\chi_1+\dotsb+a_r\chi_r);\
  a_j\ge 0\}\] is a closed convex additive cone in $i\lie a^*$.
\end{example}

Let $\mathfrak{Y}(A)$ denote the set of one-parameter subgroups $\tau
\colon \R\to A$.  In the following we identify $\mathfrak{Y}(A)$ with
the set $\Hom(\R,\lie a)$ of $\R$ linear maps from $\R$ into $\lie a$.

\begin{lemma}\label{lemmahilbert}
  Let $W$ be a completely 1-reducible $A$-module. Then
  \[
  \mathcal{S}_A(\{0\})=\bigcup_{\tau\in\mathfrak{Y}(A)}
  \mathcal{S}_{\tau(\R)}(\{0\}).
  \]
\end{lemma}

\begin{proof}
  We have to show that $\mathcal{S}_A(\{0\})\subset
  \bigcup_{\tau\in\mathfrak{Y}(A)} \mathcal{S}_{\tau(\R)}(\{0\})$.  
Let $z\in \mathcal{S}_A(0)$ and replace $W$ by the
smallest
  $A$-submodule which contains $z$.  We have the isotypic
  decomposition $W=\oplus W_{\chi_j}$ where the $\chi_j\colon \lie
  a\to \R$, $j=1,\dotsc,r$, are linear functions determining the
  weight spaces $W_{\chi_j}:= \{w\in W;\ \exp \xi\cdot
  w=\exp(\chi_j(\xi))w \text{ for all $\xi\in \lie a$}\}$. Since $W$
  is spanned by $A\cdot z$, no $\chi_j$ is identically zero, and we
  have hyperplanes $H_j:=\{\xi\in\lie a;\ \chi_j(\xi)= 0\}$.
  
  Choose a connected component $\lie a^+$ of $\lie a\setminus
  (\cup_{j=1}^r H_j)$ such that $\chi_1>0$ on $\lie a^+$.  For a
  weight $\chi_j$ we write $\chi_j>0$ if $\chi_j\vert\lie a^+>0$, and
  we write $\chi_j<0$ if $-\chi_j>0$.  We have a decomposition
  $W=W^-\oplus W^+$ where $W^-=\oplus_{\chi_j<0}W_{\chi_j}$ and
  $W^+=\oplus_{\chi_j>0}W_{\chi_j}$.  We claim that $W^-=\{0\}$.  If
  this is not the case, then $C\cap -C\not =\{0\}$ where $C$ is the
  cone generated by the $\chi_j$.  Consequently, there is a nonzero
  element $c$ in the intersection, and we can write
  $$
  c=\sum_j s_j\chi_j=-\sum_j t_j\chi_j,\quad s_j,\ t_j\in\R^+.
  $$
  For $w=w_1+\dotsb +w_r\in W_{\chi_1}+\dotsb +W_{\chi_r}$ set
  $f(w):=\prod_j \lVert w_j\rVert^{s_j+t_j}$. Then $f$ is
  $A$-invariant and continuous with $f(z)\neq 0$ and $f(0)=0$. This
  contradicts the fact that $0\in \overline{A\cdot z}$.  Thus $W=W^+$
  and for every $\xi\in \lie a^+$ we have $\lim_{t\to
    -\infty}\exp(t\xi)\cdot z=0$.  Hence $z\in
  \bigcup_{\tau\in\mathfrak{Y}(A)} \mathcal{S}_{\tau(\R)}(Y)$.
\end{proof}

\begin{corollary} \label{null cone of A} The saturation $\mathcal{S}_A(\{0\})$ is a finite union of linear subspaces of $W$.

\end{corollary}

Now let $Z$ be a holomorphic $U\c$-space with $U$-invariant K\"ahler
structure and moment map $\mu\colon Z\to \lieu^*$. Let $G$ be a closed
compatible subgroup of $U\c$ and $X$ a $G$-stable closed subset of
$Z$.

\begin{theorem} \label{HilbertMumford}  Let $\lie a \subset \liep$
  be a maximal commutative subalgebra and $A= \exp(\lie a)$ the
  corresponding subgroup of $G$. Assume that the quotient $\pi\colon
  X\to \quot X A$ exists and is regular, that $\quot XA$ parametrizes
  the closed $A$-orbits and that $\pi$ maps closed $A$-stable sets to
  closed sets.  Let $z\in X$ and let $Y\subset\overline{G\cdot z}$ be
  closed and $G$-stable. Then there is a $k\in K$, $y\in Y$ and a one
  parameter subgroup $\tau\colon\R\to A_{ y}$ such that
  $\lim_{t\to-\infty}\tau(t)k\cdot z=y$.
\end{theorem}

\begin{proof} By Proposition~\ref{propositionrichardson} we have a $k\in K$ and
  a $y\in Y$ such that $\overline{Ak\cdot z}\cap G\cdot y\neq
  \emptyset$. We may assume that the intersection contains $y$. 
Then we
  may  replace $U\c$ by the Zariski closure of $A$ in $U\c$, i.e.,
  we may assume that $U\c$ is commutative. It follows that we may change $\mu$ by
  a constant such that $y\in\m$. Hence there is an open $U\c$-stable
  neighborhood $Z_0$ of $y$ such that the analytic Hilbert quotient
  $\pi\colon Z_0\to \quot {Z_0} U\c$ exists and such that $U\c\cdot y$
  is closed in $Z_0$.  Then $A\cdot y$ is closed in $U\c\cdot y$ and
  hence in $Z_0$.  We have an $A$-slice at $y$ for the $A$-action on
  $Z_0$, hence for the $A$-action on $X$.  Thus we may assume that
  $X=\twist A {A_y}S$ where $S$ is a closed $A_y$-stable subset of an
  open $A_y$-saturated subset of $T_yZ$ relative to the quotient
  $T_yZ\to \quot {T_yZ} A_y$. Since $A_y\subset\exp i(\lieu_y)\c$ is a compatible subgroup of $(U_y)\c$, the
representation on
$T_zZ$ is completely 1-reducible, and using 
Lemma~\ref{lemmahilbert} we
  find the desired one-parameter subgroup $\tau\colon\R\to A_y$.
\end{proof}

\begin{corollary}\label{hilbert for analytic quotient} 
  Assume that the analytic Hilbert quotient $\quot Z{U\c}$ exists. Let
  $z\in Z$ and let $Y\subset \overline{G\cdot z}$ be closed and
  $G$-stable.  Then there is a $k\in K$, $y\in Y$ and a one parameter
  subgroup $\tau\colon\R\to A_{ y}$ such that
  $\lim_{t\to-\infty}\tau(t)k\cdot z=y$.
\end{corollary}
 \begin{proof} By the results in 
   section~\ref{quotients by reductive groups} we may assume that
   $Z=\mathcal{S}_{U\c}(\m)$.  Then the quotient $\pi\colon Z\to\quot
   ZA$ exists and has the desired properties. Now apply
   Theorem~\ref{HilbertMumford}
\end{proof}


\newcommand{\noopsort}[1]{} \newcommand{\printfirst}[2]{#1}
\newcommand{\singleletter}[1]{#1} \newcommand{\switchargs}[2]{#2#1}
\providecommand{\bysame}{\leavevmode\hbox to3em{\hrulefill}\thinspace}


\begin{thebibliography}{HeHuKu93}
  
\bibitem[Ab75]{Abels} H.~Abels, \emph{Parallelizability of proper
    actions, global {$K$}-slices and maximal compact subgroups}, Math.
  Ann. \textbf{212} (1975), 1--19.
  
\bibitem[AzLo92]{AzadLoeb3} H.~Azad and J.~J.~Loeb,
  \emph{Plurisubharmonic functions and K\"ahlerian metrics on
    complexification of symmetric spaces}, Indag. Math. N.S.
  \textbf{3:4} (1992), 365--375.
  
\bibitem[AzLo93]{AzadLoeb1} H.~Azad and J.~J.~Loeb,
  \emph{Plurisubharmonic functions and the Kemp-Ness Theorem}, Bull.
  London Math. Soc.  \textbf{25} (1993), 162--168.
  
\bibitem[AzLo99]{AzadLoeb2} H.~Azad and J.~J.~Loeb, \emph{Some
    applications of plurisubharmonic functions to orbits of real
    reductive groups}, Indag. Math. N.S. \textbf{10} (1999), 473--482.
  
\bibitem[Ch46]{Chevalley} C.~Chevalley, \emph{Theory of Lie Groups},
  Princeton University Press, Princeton, 1946.
  
\bibitem[GuSt91] {GuSt} V.~Guillemin and M.~Stenzel, \emph{Grauert
    tubes and the homogeneous Monge-Amp\`ere equation I}, J. Diff.
  Geom., \textbf{34} (1991), 561--570.
 
\bibitem[He91]{Heinzner} P.~Heinzner, \emph{Geometric invariant theory
    on Stein spaces}, Math. Ann. \textbf{289} (1991), 631-662.

  
\bibitem[He93]{HeinznerBull} P.~Heinzner, \emph{Equivariant
    holomorphic extensions of real analytic manifolds}, Bull. Soc.
  Math. France \textbf{121} (1993), 445--463.
  
\bibitem[HeHu96]{HeinznerHuckleberryInv} P.~Heinzner, A.~Huckleberry
  \emph{K\"ahlerian potentials and convexity properties of the moment
    map}, Invent.\ math.\ \textbf{126} (1996), 65--84.
  
\bibitem[HeHu04]{HeinznerHuckleberry} P.~Heinzner, A.~Huckleberry
  \emph{Complex geometry of Hamiltonian actions}, book, to appear.
  
\bibitem[HeHuKu93]{HeinznerHuckleberryKutschebauch}P.~Heinzner,
  A.~Huckleberry, F.~Kutzschebauch, \emph{A real analytic version of
    {A}bels' theorem and complexifications of proper {L}ie group
    actions}, in Complex Analysis and Geometry (Trento, 1993), Lecture
  Notes in Pure and Appl. Math. \textbf{173}, Dekker, New York, 1996,
  {229--273}.

  
  
\bibitem[HeHuLo94]{HeinznerHuckleberryLoose} P.~Heinzner,
  A.~Huckleberry and F.~Loose, \emph{K\"ahlerian extensions of the
    symplectic reduction}, J. reine und angew. Math. \textbf{455}
  (1994), 123--140.

  
\bibitem[HeLo94]{HeinznerLoose} P.~Heinzner and F.~Loose,
  \emph{Reduction of complex Hamiltonian $G$-spaces}, Geometric and
  Functional Analysis \textbf{4} (1994), 288--297.
  
\bibitem[HeMiPo98]{HeinznerMiglioriniPolito}P.~Heinzner, L.~Migliorini
  and M.~Polito, \emph{Semistable quotients}, Ann.  Scuola Norm. Sup.
  Pisa \textbf{26} (1998), 233--248.
  
\bibitem[Helg78]{Helgason} S.~Helgason, \emph{Differential Geometry,
    Lie Groups, and Symmetric Spaces}, Academic Press, New York, 1978.
  
\bibitem[Ho65]{Hochschild} G. Hochschild, \emph{The structure of Lie
    groups}, Holden-Day, San Francisco, 1965.

  
\bibitem[KeNe78]{KempfNess} G.~Kempf and L.~Ness, \emph{The length of
    vectors in representation spaces}, Lecture Notes in Math
  \textbf{732} (1978), Springer, New York, 233--243.
  

\bibitem[LeSz91]{LempSz} L.~Lempert and R.~Sz\"oke, \emph{Global
    solutions of the homogeneous complex Monge-Amp\`ere equation and
    complex structures on the tangent bundle of Riemannian manifolds},
  Math. Ann. \textbf{290} (1991), 689--712.
\bibitem[Lu75]{LunaReal} D.~Luna, \emph{Sur certaines op\'erations
    diff\'erentiables des groupes de {L}ie}, Amer. J. Math.
  \textbf{97} (1975), 172--181.
  
\bibitem[Mos55a]{Mostow1} G.~D.~Mostow, \emph{Some new decomposition
    theorems for semisimple groups}, Memoirs Amer. Math. Soc.
  \textbf{14} (1955), 31--54.
  
\bibitem[Mos55b]{Mostow2} G.~D.~Mostow, \emph{On covariant fiberings
    of Klein spaces}, Amer. J. Math.  \textbf{77} (1955), 247--278.

\bibitem[Nara62]{Narasimhan} R.~Narasimhan,  \emph{The Levi problem for complex spaces, II},
Math. Ann. \textbf{146} (1962),  195--216.
  
\bibitem[O'SSj00]{OsheaSjamaar} L.~O'Shea and R.~Sjamaar, \emph{Moment
    maps and Riemannian symmetric pairs}, Math. Ann. \textbf{317}
  (2000), 415--457.

  
\bibitem[Pa73]{Palais} R.~S.~Palais, \emph{On the existence of slices
    for actions of non-compact Lie groups}, Ann.\ of Math.  \textbf{73}
  (1961), 295--323.
  
\bibitem[RiSl90]{RichardsonSlodowy} R.~W.~Richardson and P.~Slodowy,
  \emph{Minimum vectors for real reductive algebraic groups}, J.\ 
  London Math.\ Soc.\ \textbf{42} (1990), 409--429.
\bibitem[Sc88]{Schwarz} G.~W. Schwarz, \emph{The topology of algebraic
    quotients}, Topological methods in algebraic transformation
  groups, edited by H. Kraft et al., Progress in Math. \textbf {80},
  Birkh\"auser Verlag, Basel-Boston, 1989, 135--152.

\end{thebibliography}
\end{document}